\documentclass[11pt]{article}
\usepackage{amsmath,amssymb}
\usepackage{graphicx}
\usepackage{subfigure}

\def\sgn{{\hbox{sgn}}}
\def\Re{{\hbox{Re}}}
\def\tr{{\hbox{\rm Tr}}}
 
\def\supp{{\hbox{\rm supp}}}
\def\C{{\hbox{\bf C}}}
\def\R{{{\mathbf{R}}}}
\def\E{{\hbox{\bf E}}}
\def\P{{\hbox{\bf P}}}

\def\eps{\varepsilon}
\newenvironment{proof}{\noindent {\bf Proof} }{\endprf\par}
\def \endprf{\hfill {\vrule height6pt width6pt depth0pt}\medskip}
\def\emph#1{{\it #1}} \def\textbf#1{{\bf #1}}

\def\RR{\mathbb{R}}

\newcommand{\bomega}{{\mbox{\boldmath{$\omega$}}}}
\newcommand{\bt}{\mathbf t}

\newtheorem{theorem}{Theorem}[section]
\newtheorem{lemma}[theorem]{Lemma}
\newtheorem{corollary}[theorem]{Corollary}

\numberwithin{equation}{section}

\def\ZZ{\mathbb{Z}}
\def\supf{T}

\begin{document}

\title{Robust Uncertainty Principles:\\
  Exact Signal Reconstruction from Highly Incomplete Frequency
  Information}



\author{Emmanuel Candes$^{\dagger}$, Justin Romberg$^{\dagger}$, and Terence Tao$^{\sharp}$\\
  \vspace{-.3cm}\\
  $\dagger$ Applied and Computational Mathematics, Caltech, Pasadena, CA 91125\\
  \vspace{-.5cm}\\
  $\sharp$ Department of Mathematics, University of California, Los
  Angeles, CA 90095
} 

\date{June 10, 2004} 

\maketitle

\begin{abstract}
  This paper considers the model problem of reconstructing an object
  from incomplete frequency samples. Consider a discrete-time signal
  $f \in \C^N$ and a randomly chosen set of frequencies $\Omega$ of
  mean size $\tau N$. Is it possible to reconstruct $f$ from the
  partial knowledge of its Fourier coefficients on the set $\Omega$?
  
  A typical result of this paper is as follows: for each $M > 0$,
  suppose that $f$ obeys
  $$
  \# \{t, \, \, f(t) \neq 0 \} \le \alpha(M) \cdot (\log N)^{-1}
  \cdot \# \Omega,
$$
then with probability at least $1-O(N^{-M})$, $f$ can be
reconstructed exactly as the solution to the $\ell_1$ minimization
problem
$$
\min_g \sum_{t = 0}^{N-1} |g(t)|, \quad \text{s.t. } \, \hat g(\omega)
= \hat f(\omega) \text{ for all } \,\, \omega \in \Omega.
$$
In short, exact recovery may be obtained by solving a convex
optimization problem.  We give numerical values for $\alpha$ which
depends on the desired probability of success; except for the
logarithmic factor, the condition on the size of the support is sharp.

The methodology extends to a variety of other setups and higher
dimensions. For example, we show how one can reconstruct a piecewise
constant (one or two-dimensional) object from incomplete frequency
samples---provided that the number of jumps (discontinuities) obeys
the condition above---by minimizing other convex functionals such as
the total-variation of $f$.

\end{abstract}


{\bf Acknowledgments.} E.~C. is partially supported by a National
Science Foundation grant DMS 01-40698 (FRG) and by an Alfred P.  Sloan
Fellowship.  J.~R.~is supported by National Science Foundation grants
DMS 01-40698 and ITR ACI-0204932.  T.~T.~is a Clay Prize Fellow and is
supported in part by grants from the Packard Foundation.  E.~C.~and
T.T.~thank the Institute for Pure and Applied Mathematics at UCLA for
their warm hospitality. E.~C.~would like to thank Amos Ron and David
Donoho for stimulating conversations, and Po-Shen Loh for early
numerical experiments on a related project.

\pagebreak

\section{Introduction}
\label{sec:introduction}


In many applications of practical interest, we often wish to
reconstruct an object (a discrete signal, a discrete image, etc.) from
incomplete Fourier samples.  In a discrete setting, we may pose the
problem as follows; let $\hat{f}$ be the Fourier transform of a
discrete object $f(t)$, $t \in\ZZ^d_N := \{0, 1,
\ldots, N-1\}^d$,
\[
\hat{f}(\omega) = \sum_{t\in\ZZ^d_N} f(t) e^{-i\omega \cdot t}.
\]
The problem is then to recover $f$ from partial frequency information,
namely, from $\hat{f}(\omega)$, where $\omega =
(\omega_1,\ldots,\omega_d)$ belongs to some set $\Omega$ of
cardinality less than $N^d$---the size of the discrete object.

In this paper, we show that we can recover $f$ {\em exactly} from
observations $\hat{f}|_\Omega$ on small set of frequencies provided
that $f$ is {\em sparse}.  The recovery consists of solving a
straightforward optimization problem that finds $f^\sharp$ of minimal
complexity with $\hat{f}^\sharp(\mathbf{\omega}) =
\hat{f}(\mathbf{\omega})$, $\forall \mathbf{\omega}\in\Omega$.


\subsection{A puzzling numerical experiment}
\label{sec:puzphantom}

\begin{figure}
\centerline{
\begin{tabular}{cc}
\includegraphics[width=2.5in]{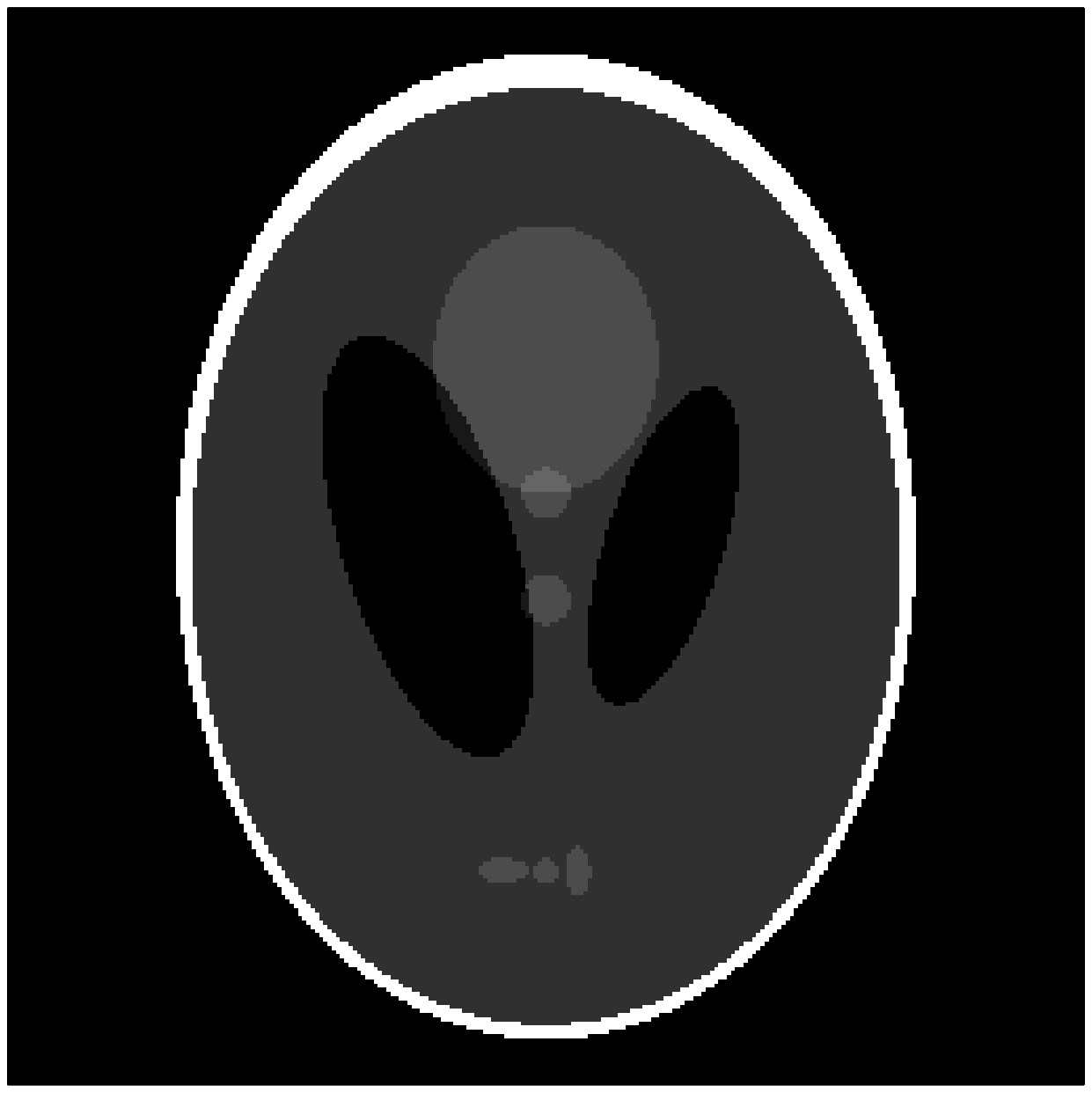} &
\includegraphics[width=2.5in]{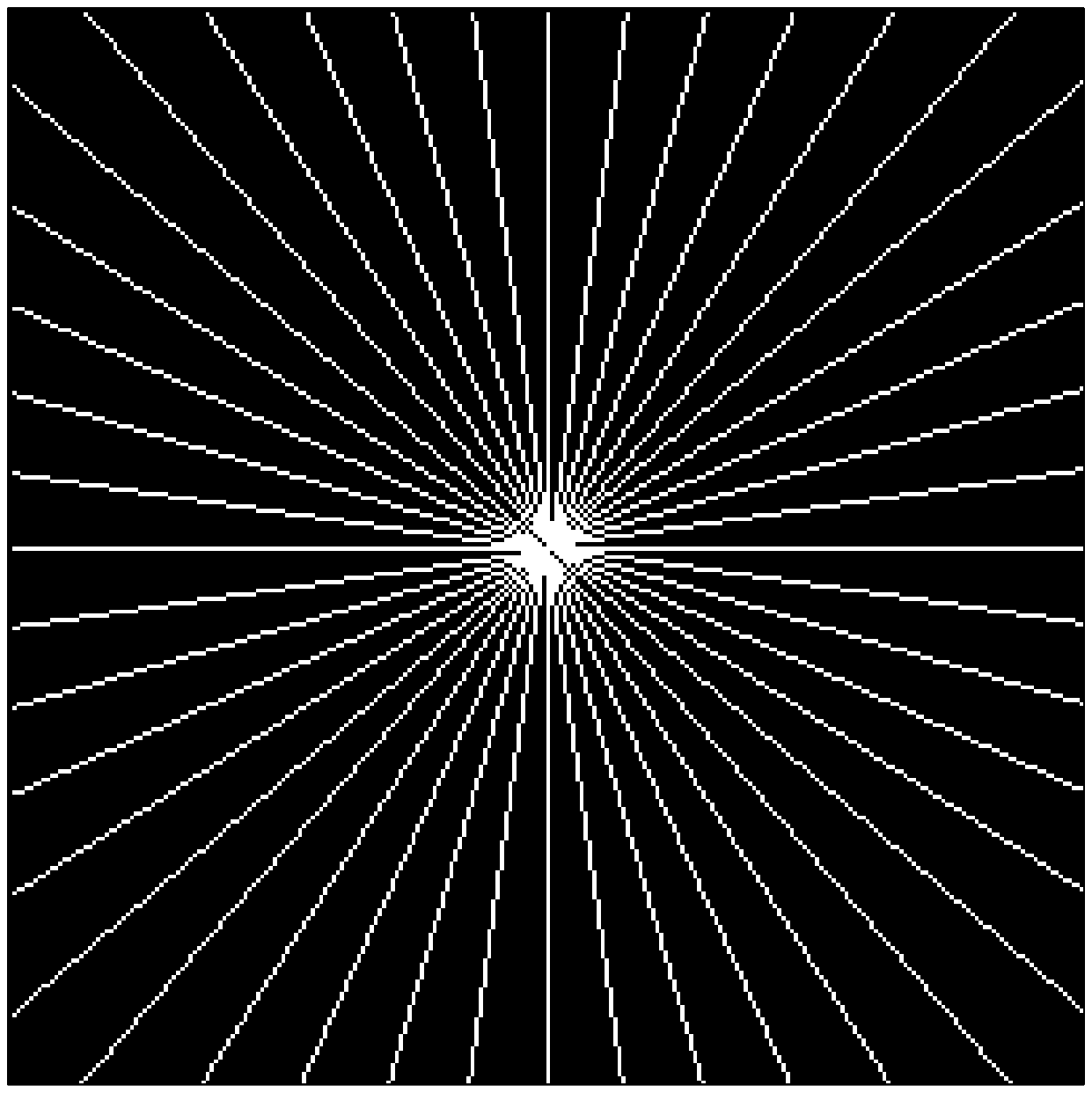} \\
(a) & (b) \\
\includegraphics[width=2.5in]{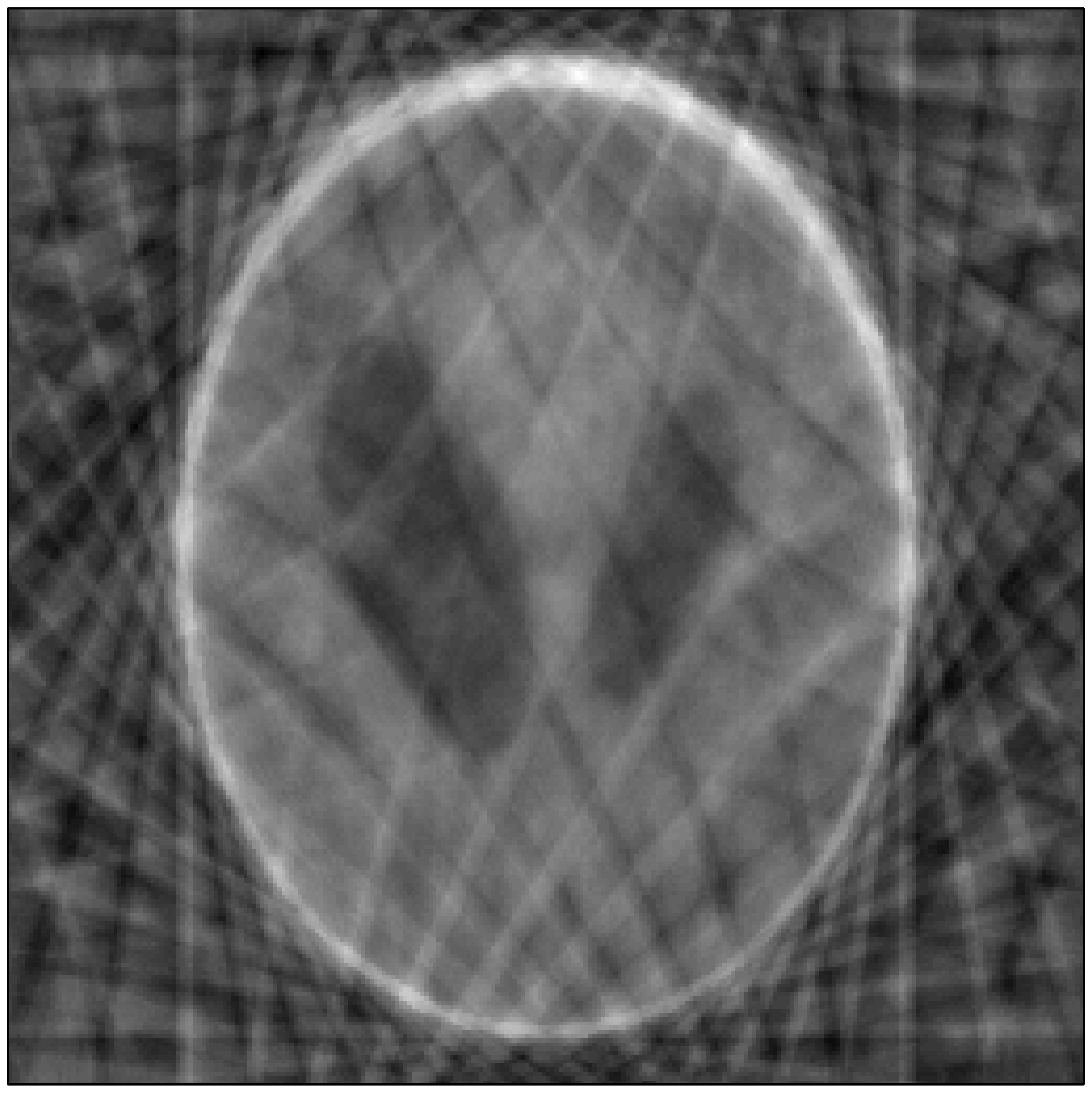} &
\includegraphics[width=2.5in]{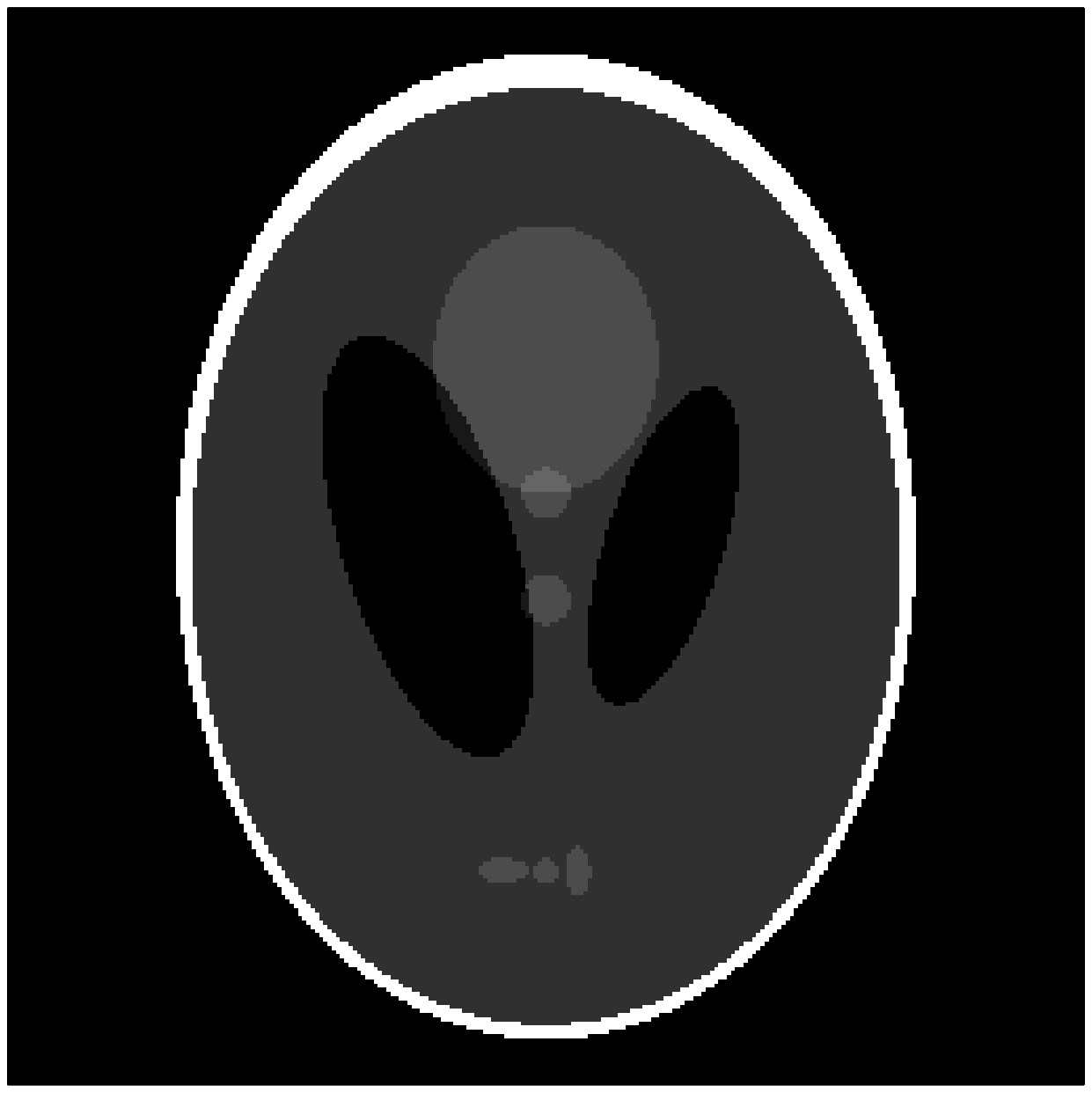} \\
(c) & (d)
\end{tabular}
}
\caption{Example of a simple recovery problem.  (a) The Logan-Shepp
  phantom test image.  (b) Sampling 'domain' in the frequency plane;
  Fourier coefficients are sampled along 22 approximately radial
  lines.  (c) Minimum energy reconstruction obtained by setting
  unobserved Fourier coefficients to zero.  (d) Reconstruction
  obtained by minimizing the total-variation, as in \eqref{eq:TV}.
  The reconstruction is an exact replica of the image in (a).}
\label{fig:phantom}
\end{figure}

This idea is best motivated by an experiment with surprisingly
positive results.
Consider a simplified version of the classical 'tomography' problem in
medical imaging: we wish to reconstruct a 2D image $f(t_1,t_2)$ from
samples $\hat{f}|_\Omega$ of its discrete Fourier transform on a
star-shaped domain $\Omega$ \cite{BreslerDelaney}.
Our choice of domain is not contrived; many real imaging devices can
collect high-resolution samples along radial lines at relatively few
angles.  
Figure~\ref{fig:phantom}(b) illustrates a typical case where one
gathers $512$ samples along each of $22$ radial lines. 

Frequently discussed approaches in the literature of medical imaging
for reconstructing an object from 'polar' frequency samples are the
so-called filtered backprojection algorithms.  In a nutshell, one
assumes that the Fourier coefficients at all of the unobserved
frequencies are zero (thus reconstructing the image of ``minimal
energy'' under the observation constraints). This strategy does not
perform very well, and could hardly be used for medical diagnostic
\cite{Mistretta}.  The reconstructed image, shown in
Figure~\ref{fig:phantom}(c), has severe nonlocal artifacts caused by
the angular undersampling.
A good reconstruction algorithm, it seems, would have to guess the
values of the missing Fourier coefficients.  In other words, one would
need to interpolate $\hat{f}(\omega_1,\omega_2)$.  This is highly
problematic, however; predictions of Fourier coefficients from their
neighbors are very delicate, due to the global and
highly oscillatory nature of the Fourier transform.
Going back to our example, we can see the problem immediately.  To
recover frequency information near $(\omega_1,\omega_2)$, where
$\omega_1$ is near $\pm\pi$, we would need to interpolate $\hat{f}$ at
the Nyquist rate $2\pi/N$.  However, we only have samples at rate
about $\pi/22$; the sampling rate is almost $50$ times smaller than
the Nyquist rate!

We propose instead a strategy based on convex optimization.
Let $\|g\|_{BV}$ be the total-variation norm of a two-dimensional
object $g$ which for discrete data $g(t_1,t_2)$, $0 \le t_1, t_2 \le N-1$,  
takes the form 
\[
\|g\| _{BV} = \sum_{t_1, t_2} 
\sqrt{|D_1 g(t_1,t_2)|^2 + |D_2 g (t_1, t_2)|^2},  
\]
where $D_1$ is the finite difference $D_1 g = g(t_1,t_2) - g(t_1 - 1,
t_2)$ and $D_2 g = g(t_1,t_2) - g(t_1,t_2 - 1)$.  To recover $f$ from
partial Fourier samples, we find a solution $f^\sharp$ to the
optimization problem 
\begin{equation}
  \label{eq:TV}
 \min \|g\|_{BV} \qquad \text{subject to} \qquad 
\hat g(\mathbf{\omega}) = 
\hat f(\mathbf{\omega}) \text{ for all } \omega \in \Omega. 
\end{equation}
In a nutshell, given partial observation $\hat f_{|\Omega}$, we seek a
solution $f^\sharp$ with minimum complexity---here Total Variation
(TV)---and whose 'visible' coefficients match those of the unknown
object $f$.  Our hope here is to partially erase some of the artifacts
classical reconstruction methods exhibit (which tend to have large TV
norm) while maintaining fidelity to the observed data via the
constraints on the Fourier coefficients of the reconstruction.

When we use \eqref{eq:TV} for the recovery problem illustrated in
Figure~\ref{fig:phantom} (with the popular Logan-Shepp phantom as a
test image), the results are surprising.  The reconstruction is {\em
  exact}; that is, $f^\sharp = f$!  Now this numerical result is not
special to this phantom.  In fact, we performed a series of
experiments of this type and obtained perfect reconstruction on many
similar test phantoms.

\subsection{Main Results}

This paper is about a quantitative understanding of this very special
phenomenon. For which classes of signals/images can we expect perfect
reconstruction? What are the trade-offs between complexity and number
of samples? In order to answer these questions, we first develop a
fundamental mathematical understanding of a special one-dimensional
model problem; we then exhibit reconstruction strategies which are
shown to exactly reconstruct the unknown signal and can be deployed in
many related and sophisticated reconstruction setups.


For a signal $f \in \C^N$, we define the classical discrete transform
Fourier transform ${\cal F} f = \hat f: \C^N \to \C^N$ by
\begin{equation}
  \label{eq:discrete-fourier}
\hat f(k) := \sum_{t = 0}^{N-1} f(t) \, e^{-i \omega_k t}, \quad \omega_k 
= \frac{2\pi k}{N}, \, k = 0, 1, \ldots, N - 1.
\end{equation}
If we are given the value of the Fourier coefficients $\hat f(k)$ for
\emph{all} frequencies $k \in \ZZ_N$, then one can
obviously reconstruct $f$ exactly via the Fourier inversion formula
\[
 f(t) = \frac{1}{N} \sum_{k = 0}^{N-1} \hat f(k) \, e^{i \omega_k t}. 
\]
Now suppose that we are only given the Fourier coefficients $\hat
f|_\Omega$ sampled in some partial subset $\Omega \subsetneq \ZZ_N$ of
all frequencies (here and below we abuse notations and identify the
frequencies $\omega_k = 2 \pi k/N$ with the corresponding integers
whenever convenient).  Of course, this is not enough information by
itself to reconstruct $f$ exactly, since $f$ has $N$ degrees of
freedom and we are only specifying $|\Omega| < N$ of those degrees
(here and below $|\Omega|$ denotes the cardinality of $\Omega$).

Suppose, however, that we also specify that $f$ is supported on a
small (but {\em a priori} unknown) subset $\supf$ of $\ZZ_N$; that is,
we assume that $f$ can be written as a sparse superposition of spikes
\[ 
f = \sum_{t\in\supf} \alpha_t \delta_t, \qquad \delta_t(t^\prime) =
1_{\{t^\prime = t\}}.
\]
If $|\supf|$ is small enough, we can recover $f$ exactly:
\begin{theorem}\label{abstract-recover}  
  Suppose that the signal length $N$ is a prime integer. Let $\Omega$ be
  a subset of $\{0, \ldots, N-1\}$, and let $f$ be a vector supported
  on $\supf$ such that
\begin{equation}\label{supp-omega}
|\supf| \leq \frac{1}{2} |\Omega|.
\end{equation}
Then $f$ can be reconstructed uniquely from $\Omega$ and $\hat
f|_\Omega$.  Conversely, if $\Omega$ is not the set of all $N$
frequencies, then there exist distinct vectors $f, g$ such that
$|\supp(f)|, |\supp(g)| \leq \frac{1}{2}|\Omega|+1$ and such that
$\hat f|_\Omega = \hat g|_\Omega$.
\end{theorem}
\begin{proof}  
  We will need the following lemma \cite{tao:uncertainty}, from which 
  we see that with knowledge of $\supf$, we can reconstruct $f$
  uniquely (using linear algebra) from $\hat f|_\Omega$:
  
  \begin{lemma}\label{unc}(\cite{tao:uncertainty}, Corollary 1.4)  
    Let $N$ be a prime integer and $\supf, \Omega$ be subsets of
    $\ZZ_N$. Put $\ell_2(\supf)$ (resp.~$\ell_2(\Omega)$) to be the
    space of signals that are zero outside of $\supf$
    (resp.~$\Omega$).  The restricted Fourier transform ${\cal
      F}_{\supf\to\Omega}: \ell_2(\supf)\to \ell_2(\Omega)$ is defined
    as
    \[
    {\cal F}_{\supf\to\Omega} f := \hat f|_\Omega 
    \hbox{ for all } f\in \ell_2(\supf),
    \]
    If $|\supf| = |\Omega|$, then ${\cal F}_{\supf\to\Omega}$ is a
    bijection; as a consequence, we thus see that ${\cal F}_{\supf\to
      \Omega}$ is injective for $|\supf| \leq |\Omega|$ and surjective for
    $|\supf| \geq |\Omega|$.  
    Clearly, the same claims hold if the
    Fourier transform ${\cal F}$ is replaced by the inverse Fourier
    transform ${\cal F}^{-1}$.
  \end{lemma}
  
  To prove Theorem~\ref{abstract-recover}, we start with the former claim.
  Suppose for contradiction that there were two objects $f, g$ such
  that $\hat f|_\Omega = \hat g|_\Omega$ and $|\supp(f)|, |\supp(g)|
  \leq \frac{1}{2} |\Omega|$.  Then the Fourier transform of $f-g$
  vanishes on $\Omega$, and $|\supp(f-g)| \leq |\Omega|$.  By Lemma
  \ref{unc} we see that ${\cal F}_{\supp(f-g) \to \Omega}$
  is injective, and thus $f-g = 0$.  The uniqueness claim follows.
  
  Now we prove the latter claim.  Since $|\Omega| < N$, we can find
  disjoint subsets $T, S$ of $\Omega$ such that $|T|, |S| \leq
  \frac{1}{2} |\Omega| + 1$ and $|T| + |S| = |\Omega| + 1$.  Let $k_0$
  be some frequency which does not lie in $\Omega$.  Applying Lemma
  \ref{unc}, we have that ${\cal F}_{T \cup S \to \Omega \cup \{k_0\}}$
  is a bijection, and thus we can find a vector $h$ supported on $T
  \cup S$ whose Fourier transform vanishes on $\Omega$ but is non-zero
  on $k_0$; in particular, $h$ is not identically zero.  The claim now
  follows by taking $f := h|_T$ and $g := -h|_S$.
\end{proof}

Note that if $N$ is not prime, the lemma (and hence the theorem)
fails, essentially because of the presence of non-trivial subgroups of
$\ZZ_N$ with addition modulo $N$; see
\cite{DonohoStark}, \cite{tao:uncertainty} for further discussion.
However, it is plausible to think that Lemma~\ref{unc} continues to
hold for non-prime $N$ if $\supf$ and $\Omega$ are assumed to be
\emph{generic} - in particular, they are not subgroups of $\ZZ_N$, or
cosets of subgroups.  If $\supf$ and $\Omega$ are selected uniformly
at random, then it is expected that the theorem holds with probability
very close to one; one can indeed presumably quantify this statement
by adapting the arguments given above but we will not do so here.
However, we refer the reader to section \ref{sec:robust} for a rapid
presentation of informal arguments pointing out in this direction.

A refinement of the argument in Theorem~\ref{abstract-recover} shows
that for fixed sets $T$, $S$, $\Omega$ in $\ZZ_N$, the
space of vectors $f, g$ supported on $T$, $S$ such that $\hat
f|_\Omega = \hat g|_\Omega$ has dimension $|T\cup S| - |\Omega|$ when
$|T \cup S| \geq |\Omega|$, and has dimension $|T \cap S|$ otherwise.
In particular, if we let $\Sigma(N_t)$ denote those vectors whose
support has size at most $N_t$, then set of the vectors in
$\Sigma(N_t)$ which cannot be reconstructed uniquely in this class
from the Fourier coefficients sampled at $\Omega$, is contained in a
finite union of linear spaces of dimension at most $2N_t - |\Omega|$.
Since $\Sigma(N_t)$ itself is a finite union of linear spaces of
dimension $N_t$, we thus see that recovery of $f$ from $\hat
f|_\Omega$ is in principle possible \emph{generically} whenever
$|\supp(f)| = N_t < |\Omega|$; once $N_t \geq |\Omega|$, however, it
is clear from simple degrees-of-freedom arguments that unique recovery
is no longer possible.  While our methods do not quite attain this
theoretical upper bound for correct recovery, our numerical
experiements suggest that they do come within a constant factor of
this bound (see Figure \ref{fig:recover512}).

Theorem~\ref{abstract-recover} asserts that $f$ can be reconstructed
from $\hat{f}|_\Omega$ if $|\supf|\leq |\Omega|/2$ (and that this
bound is the best possible).  In principle, we can recover $f$ 
exactly by solving the combinatorial optimization problem
\begin{equation}
\label{eq:(P0)}
(P_0)\quad \quad  \min_{g \in \C^N} 
\|g\|_{\ell_0}, \quad \hat{g}|_\Omega = \hat{f}|_\Omega,
\end{equation}
where $\|g\|_{\ell_0}$ is the number of nonzero terms $\# \{t,\, g(t)
\neq 0\}$.  Solving \eqref{eq:(P0)} directly is infeasible even for
modest-sized signals.  The algorithm would let $\supf$ run over all
subsets of $\{0, \ldots, N - 1\}$ of cardinality $|\supf| \leq
\frac{1}{2} |\Omega|$ and for each $\supf$, checking whether $f$ was
in the range of ${\cal F}_{\supf\to\Omega}$ or not, and then inverting
the relevant minor of the Fourier matrix to recover $f$ once $\supf$
was determined.  It is well-known that this procedure would clearly be
very computationally expensive, however, since there are exponentially
many subsets to check; for instance, for $|\Omega| \sim N/2$, this
number scales like $4^N \cdot 3^{-3N/4}$!  As an aside comment, note
that it is not clear how to make this algorithm robust, especially
since the results in \cite{tao:uncertainty} do not provide any
effective lower bound on the determinant of the minors of the Fourier
matrix, see section \ref{sec:discussion} for a discussion of this
point. 

A more computationally efficient strategy for recovering $f$ from
$\Omega$ and $\hat f|_\Omega$ is to solve the convex problem
\begin{equation}
  \label{g-l1}
(P_1)\quad\quad  \min_{g \in \C^N} 
\|g\|_{\ell_1}: = \sum_{t \in \ZZ_N} |g(t)|, 
\quad \hat{g}|_\Omega = \hat{f}|_\Omega.
\end{equation}
The key result in this paper is that the solutions to $(P_0)$ and
$(P_1)$ are {\em equivalent} for an overwhelming percentage of the
choices for $\supf$ and $\Omega$ with
$|\supf|\leq\alpha\cdot|\Omega|/\log N$ ($\alpha>0$ is a constant): in
these cases, {\em solving the convex problem $(P_1)$ recovers $f$
  exactly}.

To establish this upper bound, we will assume that the observed 
Fourier coefficients are {\em randomly sampled}.  
To make this precise, we introduce a probability
parameter $0 < \tau < 1$, and consider the sequence $(I_k)_{1 \le k
  \le N}$ of independent Bernoulli random variables
\begin{equation}
  \label{eq:Ik}
I_k =  \begin{cases} 0 & \text{with prob. } 1-\tau,\\ 
 1 & \text{with prob. } \tau.   \end{cases} 
\end{equation}
We then define the random set of frequencies $\Omega$ as
\begin{equation}\label{omega-def}
 \Omega := \{k : I_k = 1 \}.
\end{equation}
Clearly, $|\Omega|$ follows the binomial distribution and 
\begin{equation}\label{E-omega}
 \E(|\Omega|) = \tau N. 
\end{equation}
In fact, classical large deviations arguments (or the central limit
theorem) tell us that with high probability, the size of $|\Omega|$ is
very close to $\tau N$.  Our main theorem can now be stated as
follows.

\begin{theorem}\label{main}  
  Let $f \in \C^N$ be a discrete signal and $\Omega$ be the random set
  defined in \eqref{omega-def}.  For a given accuracy parameter $M$,
  if $f$ is supported on $\supf$ and
\begin{equation}\label{supp}
 |\supf| \leq \alpha (M) \cdot (\log N)^{-1} \cdot \tau N, 
\end{equation}
then with probability at least $1-O(N^{-M})$, the minimizer to the problem
\eqref{g-l1} is unique and is equal to $f$. 
\end{theorem}
In light of \eqref{E-omega} we see that \eqref{supp} is essentially
$|\supf|\sim |\Omega|$, modulo a constant and a logarithmic factor.
Indeed, an easy modification to the second part of
Theorem~\ref{abstract-recover} shows that the condition \eqref{supp}
cannot be weakened to (for instance) $|\supp(f)| \leq
(\frac{1}{2}+\eps) \tau N$, for any $\epsilon > 0$.  The paper gives
an explicit value of $\alpha(M)$, namely, $\alpha(M) \asymp
1/[29.6(M+1)]$ although we have not pursued the question of
exactly what the optimal value might be.

In Section~\ref{sec:numexp}, we present numerical results which
suggest that in practice, we can expect to recover $f$ more than
$50\%$ of the time if $|\supf|\leq |\Omega|/4$.  For $|\supf|\leq
|\Omega|/8$, the recovery rate is above $90\%$.  Empircally, the
constants $1/4$ and $1/8$ do not seem to vary for $N$ in the range of
a few hundred to a few thousand.

\subsection{For Almost Every $\Omega$}

As the theorem suggests, there exist sets $\Omega$ and functions $f$
for which the $\ell_1$-minimization procedure does not recover $f$
correctly, even if $|\supp(f)|$ is much smaller than $|\Omega|$.  We
sketch two counter-examples:
\begin{itemize}
\item {\em Dirac's comb.} Suppose that $N$ is a perfect square and
  consider the picket-fence signal which consists of spikes of unit
  height and with uniform spacing equal to $\sqrt{N}$. This signal is
  often used as an extremal point for uncertainty principles
  \cite{DonohoStark,DonohoHuo} as one of its remarkable properties is
  its invariance through the Fourier transform. Hence suppose that
  $\Omega$ is the set of all frequencies but the multiples of
  $\sqrt{N}$, namely, $|\Omega| = N - \sqrt{N}$. Then $\hat
  f|_{\Omega} = 0$ and obviously the reconstruction is identically
  zero.
  
  Note that the problem here does not really have anything to do with
  $\ell_1$-minimization per se; $f$ cannot be reconstructed from its
  Fourier samples on $\Omega$ thereby showing that Theorem
  \ref{abstract-recover} does not work 'as is' for arbitrary sample
  sizes.
    
\item {\em Box signals}. The example above suggests that in some sense
  $|T|$ must not be greater than about $\sqrt{|\Omega|}$. In fact,
  there exist more extreme examples. Assume the sample size $N$ is
  large and consider for example the indicator function $f$ of the
  interval $T := \{ t: -N^{-0.01} < t < N^{0.01} \}$ and let $\Omega$
  be the set $\Omega := \{ k: N/3 < k < 2N/3 \}$.  Let $h$ be a
  function whose Fourier transform $\hat h$ is a non-negative bump
  function adapted to the interval $\{k : -N/6 < k < N/6 \}$
  which equals 1 when $-N/12 < k < N/12$.  Then $|h(t)|^2$ has Fourier
  transform vanishing in $\Omega$, and is rapidly decreasing away from
  $t=0$; in particular we have $|h(t)|^2 = O(N^{-100})$ for $t \not
  \in T$.  On the other hand, one easily computes that $|h(0)|^2 > c$
  for some absolute constant $c > 0$.  Because of this, the signal $f
  - \eps |h|^2$ will have smaller $\ell_1$-norm than $f$ for $\eps >
  0$ sufficiently small (and $N$ sufficiently large), while still
  having the same Fourier coefficients as $f$ on $\Omega$.  Thus in
  this case $f$ is not the minimizer to the problem $(P_1)$, despite
  the fact that the support of $f$ is much smaller than that of
  $\Omega$.
\end{itemize}

The above counterexamples relied heavily on the special choice of
$\Omega$ (and to a lesser extent of $\supp(f)$); in particular, it
needed the fact that the complement of $\Omega$ contained a large
interval (or more generally, a long arithmetic progression). But for
most sets $\Omega$, large arithmetic progressions in the complement do
not exist, and the problem largely disappears. In short, Theorem
\ref{main} essentially says is that for {\em most} sets $|\supf| \sim
|\Omega|$, the inequality holds.


\subsection{Extensions}

As mentioned earlier, results on our model problem extend easily to
higher dimensions as well as to other setups. To be concrete consider
the problem of recovering a one-dimensional piecewise constant signal
via
\begin{equation}
  \label{eq:minTV-1D}
\min_g \,\,\, \sum_{t \in \ZZ_N} |g(t) - g(t-1)| \qquad 
\hat g|_\Omega = \hat f|_\Omega, 
\end{equation}
where we adopt the convention that $g(-1) = g(N-1)$. In a nutshell,
model \eqref{g-l1} is obtained from \eqref{eq:minTV-1D} after
differentiation.  Indeed, let $\delta$ be the vector of first difference
$\delta(t) = g(t) - g(t-1)$, and note that $\sum \delta(t) = 0$.
Obviously,
\[
\hat \delta(\omega) = (1 - e^{-i\omega}) \hat g(\omega), \quad
\text{for all } \omega \neq 0  
\]
and, therefore, with $\upsilon(\omega) = (1 - e^{-i\omega})^{-1}$, the
problem is identical to
\[
\min_\delta \,\, \|\delta\|_{\ell_1} \qquad \hat \delta|_{\Omega
  \setminus \{0\}} = (\upsilon \hat f)|_{\Omega \setminus \{0\}}, \,\,
\hat \delta(0) = 0,
\]
which is precisely what we have been studying. 
\begin{corollary}
  \label{main-corollary}
  Put $T = \{t,\, f(t) \neq f(t-1)\}$. Under the assumptions of
  Theorem \ref{main}, the minimizer to the problem \eqref{eq:minTV-1D}
  is unique and is equal $f$ with probability at least
  $1-O(N^{-M})$---provided, of course, that $f$ be adjusted so that
  $\sum f(t) = \hat f(0)$.
\end{corollary}

We now explore versions of Theorem \ref{main} in higher dimensions. To
be concrete, consider the two-dimensional situation (statements in
arbitrary dimensions are exactly of the same flavor):
\begin{theorem}\label{main2}  
  Put $N = n^2$. We let $f(t_1, t_2), 1 \le t_1, t_2 \le n$ be a
  discrete signal and $\Omega$ be the random set defined as in
  \eqref{omega-def}. Assume that for a given accuracy parameter $M$,
  $f$ is supported on $\supf$ obeying \eqref{supp}. Then with
  probability at least $1-O(N^{-M})$, the minimizer to the problem
  \eqref{g-l1} is unique and is equal to $f$.
\end{theorem}
We will not prove this result as the strategy is exactly parallel to
that of Theorem \ref{main}.  Just as in the one-dimensional case, a
similar statement for piecewise constant functions exists provided, of
course, that the support of $f$ be replaced by $\{(t_1, t_2): |D_1
f(t_1,t_2)|^2 + |D_2 f (t_1, t_2)|^2 \neq 0\}$.  We omit the details.

We hope that we managed to suggest that there actually are a variety
of results similar to Theorem \ref{main}, and we only selected a few
instances. As a matter of fact, those provide a precise quantitative
understanding of the `surprising result' discussed at the beginning of
this paper.


\subsection{Relationship to Uncertainty Principles}

From a certain point of view, our results are connected to the
so-called {\em uncertainty principles} \cite{DonohoStark,DonohoHuo}
which say that it is difficult to localize a signal $f \in \C^N$ both
in time and frequency at the same time. Indeed, classical arguments
show that $f$ is the unique minimizer of $(P_1)$ if and only if
\[
\sum_{t \in \ZZ_N} |f(t) + h(t)| > \sum_{t \in \ZZ_N} |f(t)|, \quad
\forall h \neq 0, \, \hat h|_{\Omega} = 0
\]
Put $T = \supp(f)$ and apply the triangle inequality 
\[
\sum_{\ZZ_N} |f(t) + h(t)| = \sum_{T} |f(t) + h(t)| + \sum_{T^c} |h(t)| \ge
\sum_{T} |f(t)| - |h(t)| + \sum_{T^c} |h(t)|. 
\]
Hence, a sufficient condition to establish that $f$ is our unique
solution would be to show that 
\[
\sum_T |h(t)| < \sum_{T^c} |h(t)| \qquad \forall h \neq 0, \, \hat
h|_{\Omega} = 0. 
\]
or equivalently $\sum_T |h(t)| < \frac{1}{2} \|h\|_{\ell_1}$. The
connection with the uncertainty principle is now explicit; $f$ is the
unique minimizer if it is impossible to `concentrate' half of the
$\ell_1$ norm of a signal that is missing frequency components in
$\Omega$ on a 'small' set $T$. For example, \cite{DonohoStark}
guarantees exact reconstruction if
\[
2 |T| \cdot  (N - |\Omega|) < N.
\]
Take $|\Omega| < N/2$, then that condition says that $|T|$ must be
zero which, of course, is far from being the content of Theorem
\ref{main}. In truth, this paper does not follow this classical
approach.  Instead, we will use duality theory to study the solution
of $(P_1)$.

\subsection{Robust Uncertainty Principles} 
\label{sec:robust}

Underlying our analysis is a new notion of uncertainty principle which
holds for almost any pair $(\supp(f), \supp(\hat f))$. With $T =
\supp(f)$ and $\Omega = \supp(\hat f)$, the classical discrete
uncertainty principle \cite{DonohoStark} says that
\begin{equation}
  \label{eq:UP-Donoho}
  |T| + |\Omega| \ge 2 \sqrt{N}.
\end{equation}
with equality obtained for signals such as the Dirac's comb. As we
mentioned above, such extremal signals correspond to very special
pairs $(T, \Omega)$. However, for most choices of $T$ and $\Omega$,
the analysis presented in this paper shows that it is {\em impossible}
to find $f$ such that $T = \supp(f)$ and $\Omega = \supp(\hat f)$
unless
\begin{equation}
  \label{eq:Robust-UP}
  |T| + |\Omega| \ge \gamma(M)  \cdot (\log N)^{-1/2} \cdot N,   
\end{equation}
which is considerably stronger than \eqref{eq:UP-Donoho}. Here, the
statement 'most pairs' says again that the probability of selecting a
random pair $(T, \Omega)$ violating \eqref{eq:Robust-UP} is at most
$O(N^{-M})$. (We are of course aware of numerical studies in
\cite{DonohoStark} pointing out the lack of sharpness of the
uncertainty principle when $T$ is random.)

In some sense, \eqref{eq:Robust-UP} is the typical uncertainty
relation one can generally expect (as opposed to
\eqref{eq:UP-Donoho}), hence, justifying the title of this paper. 
Because of space limitation, we are unable to belaborate on this fact
and its implications any further, but will do so in a companion paper.

\subsection{Connections with existing work}

The idea of relaxing a combinatorial problem into a convex problem is
not new and goes back a long way. For example,
\cite{DobsonSantosa,SantosaSymes} used the idea of minimizing $\ell_1$
norms to recover spike trains. The motivation is that this makes
available a host of computationally feasible procedures. For example,
a convex problem of the type (\ref{g-l1}) can be practically solved
using techniques of linear programming such as interior point methods
\cite{BP}.

Now, there exists some evidence that in special situations the unique
solution to an $\ell_1$ minimization problem coincides with that of
the unique minimizer of the $\ell_0$ problem. For example, a series of
beautiful papers \cite{DonohoHuo,DonohoElad,EladBruckstein,
  FeuerNemirovsky,GribonvalNielsen} is concerned with a special setup
where one is given a dictionary $D$ of vectors (waveforms) of $\C^N$,
$D = (d_k)_{1 \le k \le M}$ and one seeks sparse representations of a
signal $f \in \C^N$ as a superposition of elements of
$D$
\begin{equation}
  \label{eq:Da}
  f = D \alpha.
\end{equation}
Suppose that the number of elements $M$ from $D$ is greater than the
sample size $N$, then there are many ways in which one can represent
$f$ as a superposition of elements from $D$ and one would want to find
the 'sparsest' one.  Consider the solution which minimizes the
$\ell_0$ norm of $\alpha$ subject to the constraint \eqref{eq:Da} and
that which minimizes the $\ell_1$ norm. A typical result of this body
of work is as follows: suppose that $s$ can be synthesized out of very
few elements from $D$, then the solution to both problems are unique
and are equal. We also refer to \cite{Tropp03,Tropp04} for very recent results
along these lines.

This literature certainly influenced our thinking in the sense it made
us suspect that results such as Theorem \ref{main} were actually
possible.  However, we would like to emphasize that the claims
presented in this paper are of a substantially different nature. We
give essentially two reasons:
\begin{itemize}
\item First, our model problem is different since we need to
  'guess' a signal from incomplete data, as opposed to finding the
  sparsest expansion of a fully specified signal.
\item And second, our approach is decidedly probabilistic---as opposed
  to deterministic---and thus calls for very different techniques. For
  example, underlying our analysis are delicate estimates about the
  size of random matrices, which may be of independent interest.
\end{itemize}

Besides the wonderful properties of $\ell_1$, there is a second line
of research connected to our findings. We can think of recovering a
sparse superposition of spikes from an incomplete set of observations
in the Fourier domain as a spectral estimation problem proviso
swapping time and frequency: $\hat{f}$ is a superposition of a few
complex sinusoids whose frequency and amplitude we need to determine
from a few samples.
From this point of view, our work is related to
\cite{BreslerFengI,BreslerFengII} and \cite{VetterliMarzBlu} where the
authors study sampling patterns allowing the exact reconstruction of a
signal.  These references show that the locations and amplitudes of a
sequence of $|T|$ spikes can be recovered exactly from $2|T|+1$
consecutive Fourier coefficients (in \cite{VetterliMarzBlu} for
example, the recovery requires solving a system of equations and
factoring a polynomial). Our results, namely,
Theorems~\ref{abstract-recover} and \ref{main} are quite distinct and
far more general since they address the radically different situation
in which we do not have the freedom to choose the samples at our
convenience.

Finally, it is interesting to note that our results and the references
above are also related to recent work \cite{GilbertStrauss} in finding
near-best $B$-term Fourier approximations (which is in some sense the
dual to our recovery problem).  The algorithm in
\cite{GilbertStrauss,GilbertStraussII}, which operates by estimating
the frequencies present in the signal from a small number of randomly
placed samples, produces with high probability an approximation in
sublinear time with error within a constant of the best $B$-term
approximation. First, in \cite{GilbertStraussII} the samples are again
selected to be equispaced whereas we are not at liberty to choose the
frequency samples at all since they are specified {\em a priori}. And
second, we wish to produce as a result an entire signal or image of
size $N$, so a sublinear algorithm is an impossibility.

\section{Strategy}

It is clear that at least one minimizer to $(P_1)$ exists.  On the
other hand, it is not apparent why this minimizer should be unique,
and why it should equal $f$.  In this section, we outline our strategy
for answering these questions.  Using duality theory, we will be able
to derive necessary and sufficient conditions for $(P_1)$ to recover
$f$.  We note that a similar duality approach was independently
developed in \cite{FuchsDual} for finding sparse approximations from
general dictionaries.

\subsection{Duality}

To get a feel for the line of argumentation, consider first the case
where $f$ is real-valued.  Then \eqref{g-l1} can be written as the
linear program
\begin{equation}
\label{eq:LP}
\min_{\substack{g^{+},g^{-}\in\RR^N\\ g^{+},g^{-}\geq 0}}~
\sum_{t=0}^{N-1} (g^{+}(t) + g^{-}(t)),  \quad 
\mathcal{F}_\Omega(g^{+}-g^{-}) = \hat{f}|_\Omega
\end{equation}
where $g^{+}(t) = \max(g(t),0)$, $g^{-}(t) = -\min(g(t),0)$, and
the matrix $\mathcal{F}_\Omega$ contains only the rows of the Fourier transform
matrix corresponding to entries in $\Omega$.  The corresponding
Lagrangian is
\begin{equation}
L(g^{+},g^{-}; \lambda, \mu^{+}, \mu^{-}) = 
\sum_{t=0}^{N-1} (g^{+}(t) + g^{-}(t)) ~+~ 
\lambda^H(\hat{f}|_\Omega - \mathcal{F}_\Omega(g^{+}-g^{-})) ~+~
\mu^{+*}g^{+} ~+~ (\mu^{-})^* g^{-}
\end{equation}
with $\mu^{+},\mu^{-}\geq 0$. 
At a minimum $(\tilde{g}^{+},\tilde{g}^{-})$, there will be a saddle
point in $L$, and we will have 
\begin{eqnarray*}
\mathcal{F}_\Omega(\tilde{g}^{+}-\tilde{g}^{-}) & = & \hat{f}|_\Omega \\
(\mu^{+})^* \tilde{g}^{+} & = & 0 \\
(\mu^{-})^* \tilde{g}^{-} & = & 0 \\
\frac{\partial L}{\partial \tilde{g}^{+}(t)} ~=~ I_{\{\tilde{g}^{+}(t)>0\}} 
~-~ \mathcal{F}^*_\Omega\lambda ~+~ \mu^{+} & = & 0 \\
\frac{\partial L}{\partial \tilde{g}^{-}_t} ~=~ I_{\{\tilde{g}^{-}(t)>0\}} ~+~ 
\mathcal{F}^*_\Omega\lambda ~+~ \mu^{-} & = & 0.
\end{eqnarray*}
Then for $f$ to be the minimum of \eqref{eq:LP}, we need
\begin{eqnarray}
(\mathcal{F}^*_\Omega\lambda)(t) & = & \sgn(f)(t) \quad\quad t\in\supf \\
1 - (\mathcal{F}^*_\Omega\lambda)(t) - \mu^{+} & = & 0 ~\quad\quad\quad\quad\quad t\in\supf^c \\
1 + (\mathcal{F}^*_\Omega\lambda)(t) - \mu^{-} & = & 0 ~\quad\quad\quad\quad\quad t\in\supf^c
\end{eqnarray}
with $\mu^{+},\mu^{-}\geq 0$.  In fact, for $f$ to be the unique
minimizer of \eqref{eq:LP}, it is necessary and sufficient for there
to exist a $\lambda$ such that for $P(t) = (\mathcal{F}^*_\Omega\lambda)(t)$, we 
have
\begin{eqnarray}
P(t) & = & \sgn(f)(t) \quad\quad t\in\supf \\
|P(t)| & < & 1 ~\quad\quad\quad\quad\quad t\not\in\supf.
\end{eqnarray}

Thus, to show that $f^\sharp$ is unique and is equal to $f$, it
suffices to find a trigonometric polynomial $P$ whose Fourier
transform is supported in $\Omega$---in other words, which only uses
frequencies in $\Omega$---and which matches $\sgn(f)$ on $\supp(f)$,
and has magnitude strictly less than 1 elsewhere.  The following lemma
generalizes for the case where $f$ is complex-valued.
\begin{lemma}\label{duality}  
  Let $\Omega \subset \ZZ_N$.  For a vector $f\in\C^N$, define the
  'sign' vector $\sgn(f)$ by $\sgn(f)(t) := f(t)/|f(t)|$ when $t \in
  \supp(f)$ and $\sgn(f) = 0$ otherwise.  Suppose there exists a
  vector $P$ whose Fourier transform $\hat P$ is supported in $\Omega$
  such that
  $$ P(t) = \sgn(f)(t) \hbox{ for all } t \in \supp(f) $$
  and
  $$
  |P(t)| < 1 \hbox{ for all } t \not \in \supp(f).  $$
  \begin{itemize}
  \item Then if ${\cal F}_{\supp(f) \to \Omega}$ is injective, the
    minimizer $f^\sharp$ to the problem $(P_1)$ \eqref{g-l1} is unique
    and is equal to $f$.
  \item Conversely, if $f$ is the unique minimizer of $(P_1)$, then
    there exists a vector $P$ with the above properties.
\end{itemize}
\end{lemma}

\begin{proof}  
  We may assume that $\Omega$ is non-empty and that $f$ is non-zero
  since the claims are trivial otherwise.
  
  Suppose first that such a function $P$ exists.  Let $g$ be any
  vector not equal to $f$ with $\hat{g}|_\Omega = \hat{f}|_\Omega$.
  Write $h := g-f$, then $\hat h$ vanishes on
  $\Omega$.  Observe that for any $t \in \supp(f)$ we have
\begin{align*}
 |g(t)| &= |f(t) + h(t)| \\
&= ||f(t)| + h(t) \, \overline{\sgn(f)(t)}| \\
&\geq |f(t)| + \Re( h(t) \, \overline{\sgn(f)(t)} )\\
&= |f(t)| + \Re( h(t) \, \overline{P(t)} )
\end{align*}
while for $t \not \in \supp(f)$ we have $|g(t)| = |h(t)| \geq
\Re( h(t) \overline{P(t)} )$ since $|P(t)| < 1$.  Thus
$$
\| g \|_{\ell_1} \geq \| f \|_{\ell_1} + \sum_{t = 0}^{N-1} \Re(
h(t) \, \overline{P(t)} ).$$
However, the Parseval's formula gives
$$
\sum_{t = 0}^{N-1} \Re( h(t) \overline{P(t)} ) = \frac{1}{N}
\sum_{k = 0}^{N-1} \Re( \hat h(k) \, \overline{\hat P(k)}) = 0$$
since
$\hat P$ is supported on $\Omega$ and $\hat h$ vanishes on $\Omega$.
Thus $\|g\|_{\ell_1} \geq \|f\|_{\ell_1}$.  Now we check when equality
can hold, i.e. when $\|g\|_{\ell_1} = \| f \|_{\ell_1}$.  An
inspection of the above argument shows that this forces $|h(t)| =
\Re(h(t) \overline{P(t)})$ for all $t \not \in \supp(f)$.  Since
$|P(t)| < 1$, this forces $h$ to vanish outside of $\supp(f)$.  Since
$\hat h$ vanishes on $\Omega$, we thus see that $h$ must vanish
identically (this follows from the assumption about the injectivity of
${\cal F}_{\supp(f) \to \Omega}$) and so $g=f$.  This shows that $f$
is the unique minimizer $f^\sharp$ to the problem \eqref{g-l1}.

Conversely, suppose that $f = f^\sharp$ is the unique minimizer to
\eqref{g-l1}.  Without loss of generality we may
normalize $\|f\|_{\ell_1} = 1$.  Then the closed unit ball $B := \{ g:
\|g\|_{\ell_1} \leq 1 \}$ and the affine space $V := \{ g: \hat
g|_\Omega = \hat f|_\Omega \}$ intersect at exactly one point, namely
$f$.  By the Hahn-Banach theorem we can thus find a function $P$ such
that the hyperplane $\Gamma_1 := \{ g: \sum \Re(g(t) \, \overline{P(t)}) = 1
\}$ contains $V$, and such that the half-space $\Gamma_{\leq 1} := \{
g: \sum \Re(g(t) \, \overline{P(t)}) \leq 1 \}$ contains $B$.  By perturbing
the hyperplane if necessary (and using the uniqueness of the
intersection of $B$ with $V$) we may assume that $\Gamma_1 \cap B$ is
contained in the minimal facet of $B$ which contains $f$, namely $\{ g
\in B: \supp(g) \subseteq \supp(f)\}$.

Since $B$ lies in $\Gamma_{\leq 1}$, we see that $\sup_t |P(t)| \leq
1$; since $f \in \Gamma_1 \cap B$, we have $P(t) = \sgn(f)(t)$ when $t
\in \supp(f)$.  Since $\Gamma_1 \cap B$ is contained in the minimal
facet of $B$ containing $f$, we see that $|P(t)| < 1$ when $t \not \in
\supp(f)$.  Since $\Gamma_1$ contains $V$, we see from Parseval that
$\hat P$ is supported in $\Omega$.  The claim follows.
\end{proof}

Since the space of functions with Fourier transform supported in
$\Omega$ has $|\Omega|$ degrees of freedom, and the condition that $P$
match $\sgn(f)$ on $\supp(f)$ requires $|\supp(f)|$ degrees of
freedom, one now expects heuristically (if one ignores the open
conditions that $P$ has magnitude strictly less than 1 outside of
$\supp(f)$) that $f^\sharp$ should be unique and be equal to $f$ whenever
$|\supp(f)| \ll |\Omega|$; in particular this gives an explicit
procedure for recovering $f$ from $\Omega$ and $\hat f|_\Omega$.

\subsection{Architecture of the Argument}

Equipped with our duality theorem, we are now in a position to present
the main ideas of the argument. Fix $f$.  We may assume that $\tau N >
M \log N$ since the claim is vacuous otherwise (as we will see,
$\alpha(M) = O(1/M)$ and thus \eqref{supp} will force $f \equiv 0$, at
which point it is clear that the solution to $(P_1)$ is equal to $f =
0$).

We let $T \subset \ZZ_N$ denote the support of $f$, $T := \supp(f)$.
Let $\Omega$ be the random set defined by \eqref{omega-def}.  Since
$\tau N > M \log N$, a typical application of the large deviation
theorem shows that the cardinality of $\Omega$ is if course close to
that of its expected value, e.g. 
\begin{equation}
  \label{eq:tail-binomial}
\P(|\Omega| < \E |\Omega| - t) \le \exp(-t^2/2 \E |\Omega|).
\end{equation}
Slightly more precise estimates are possible, see \cite{MassartSharp}.
It then follows that 
\begin{equation}
  \label{eq:tail-binomial2}
\P(|\Omega| < (1-\epsilon_M) |\tau N|) \le N^{-M}, \qquad 
\epsilon_M := \sqrt{\frac{2M \log N}{|\tau N|}}. 
\end{equation}
In the sequel it will be convenient to denote by $B_M$ the event
$\{|\Omega| < (1-\epsilon_M) |\tau N|\}$.

In light of Lemma \ref{duality}, it suffices ---with probability $1 -
O(N^{-M})$--- to (1) show that the matrix ${\cal F}_{\supp(f) \to
  \Omega}$ has full rank, and (2) construct a trigonometric
polynomial $P(t)$, $0 \le t \le N -1$, whose Fourier transform is
supported on $\Omega$, matches $\sgn(f)$ on $T$, and has magnitude
strictly less than 1 outside of $T$.  To do this we shall need some
auxiliary linear transformations (i.e. matrices) as we will see next.

In this section, we will work with vectors restricted to the set $T$
and it will be convenient to let $\ell_2(T)$ denote the subspace of
such restrictions (and similarly $\ell_2(\ZZ_N) := \C^N$).
With these notations, we let $H: \ell^2(T) \to \ell_2(\ZZ_N)$ denote
the linear transform defined by
\begin{equation}\label{T-def}
Hf(t) := -\sum_{\omega \in \Omega} \,\,\, \sum_{t' \in T: t' \neq t} 
e^{i \omega (t - t')} \, f(t').
\end{equation}
Let $\iota: \ell^2(T) \to \ell_2(\ZZ_N)$ be the obvious embedding of
$\ell^2(T)$ into $\ell_2(\ZZ_N)$ (extending by zero outside of $T$),
and let $\iota^*: \ell_2(\ZZ_N) \to \ell^2(T)$ be the dual restriction
map, thus $\iota^* f := f|_T$.  Observe that $\iota^* \iota: \ell^2(T)
\to \ell^2(T)$ is simply the identity operator on $\ell^2(T)$, and
that the operator $\iota^* H: \ell^2(T) \to \ell^2(T)$ is
self-adjoint.

The key point is that the terms in \eqref{T-def} are rather
oscillatory, since we have stripped out the non-oscillatory diagonal
$t = t'$; indeed, the main idea of the argument will be to use the
randomization of $\Omega$ to treat $H$ as a ``white noise'' operator
whose eventual effect will be negligible, especially if $H$ is raised
to a high power.

To see the relevance of the operator $H$ to our problem, observe that
for all $f \in \ell^2(T)$ 
$$
(\iota - \frac{1}{|\Omega|} H)f(t) = \frac{1}{|\Omega|}
\sum_{\omega \in \Omega} \sum_{t' \in T} e^{i \omega(t - t')} f(t') =
  \frac{1}{|\Omega|} \sum_{\omega \in \Omega} \hat f(\omega) \,
  e^{i\omega t}, 
  $$
  with $\hat f(\omega)$ the Fourier coefficient of $f$ evaluated at
  the frequency $\omega$.  In particular, $(\iota - \frac{1}{|\Omega|}
  H)f$ has Fourier transform supported in $\Omega$.  Next, suppose for
  the moment that the self-adjoint operator $\iota^* \iota -
  \frac{1}{|\Omega|} \iota^* H$ from $\ell^2(T)$ to itself is
  invertible, and then set $P(t)$, $0 \le t \le N-1$, to be the
  trigonometric polynomial
\begin{equation}\label{P-def}
 P := (\iota - \frac{1}{|\Omega|}H) 
(\iota^* \iota - \frac{1}{|\Omega|}\iota^* H)^{-1} \iota^* \sgn(f).
\end{equation}

Then by the preceding discussion:
\begin{itemize}
\item{\em Frequency support.} $P$ has Fourier transform supported in
  $\Omega$;
\item{\em Spatial interpolation.} $P$ obeys
$$
\iota^* P = (\iota^* \iota - \frac{1}{|\Omega|}\iota^* H) (\iota^*
\iota - \frac{1}{|\Omega|}\iota^* H)^{-1} \iota^* \sgn(f) = \iota^*
\sgn(f),
$$
and so $P$ agrees with $\sgn(f)$ on $T$.
\end{itemize}

Consider now the invertibility issue. By definition 
$$
\iota^* \iota - \frac{1}{|\Omega|} \iota^* H = \frac{1}{|\Omega|}
[{\cal F}_{T \to \Omega}]^* {\cal F}_{T \to \Omega}. 
$$
Hence, the invertibility of $\iota^* \iota - \frac{1}{|\Omega|}
\iota^* H$ implies that ${\cal F}_{T \to \Omega}$ be injective. In
summary, to prove the theorem it will suffice to show that:
\begin{itemize}
\item{\em  Invertibility.} The operator $\iota^* \iota -
  \frac{1}{|\Omega|} \iota^* H$ is invertible (with probability $1 -
  O(N^{-M})$).
\item{\em  Magnitude on $T^c$.} The function $P$ defined in
  \eqref{P-def} obeys the bound $\sup_{t \in T^c} |P(t)| < 1$ (with
  probability $1 - O(N^{-M})$).
\end{itemize}
We first consider the former claim.


\section{Construction of the Dual Polynomial}

\subsection{Invertibility}

We would like to establish invertibility of the matrix $\iota^* \iota
- \frac{1}{|\Omega|} \iota^* H$ with high probability.  One obvious
way to proceed would be to show that the operator norm or equivalently
the largest eigenvalue of $\iota^* H$ is less than $|\Omega|$.  This
is easily done if $|\supp(f)|$ is extremely small (e.g.  much less
than $\sqrt{|\Omega|}$), simply by estimating the operator norm
directly by the Frobenius norm $\|\cdot\|_F$, which is easy to
compute explicitly. Recall that for any squared matrix $M$, the
\emph{Frobenius} norm $\|M\|_{F}$ of $M$ is defined by the formula
$$
\|M\|^2_{F} := \tr(M M^*) = \sum_{i,j} |M(i,j)|^2, 
$$
and obeys $\|M\| \le \|M\|_F$.  However, this simple approach does
not work well when $|\supp(f)|$ is large, say equal to $\alpha \cdot
(\log N)^{-1} \cdot |\Omega|$.  In this case, we have to resort to
estimating the Frobenius norm of a \emph{large power} of $\iota^* H$,
taking advantage of cancellations arising from the randomness of the
matrix coefficients of $\iota^* H$.

We state the key estimate of this section. 
\begin{theorem}
\label{moment-bound}
Put $H_0 = \iota^* H$ for short, where $H$ is the operator defined by
\eqref{T-def}. Set $c_\tau := e \log((1-\tau)/\tau)$ and let
$$
a_n = (2n-1)^{2n} \, c_\tau^{-(2n-1)}\, N \, |T|^{2n}, \quad b_n =
\frac{(2n)!}{n! \, 2^n} \, \left(\frac{\tau}{1-\tau}\right)^n \, N^n
\, |T|^{n+1}.
$$
Then
  \begin{equation}
    \label{eq:moment-bound}
\E[\tr(H_0^{2n})] \le n \, \left(\frac{1+\sqrt{5}}{2}\right)^{2n} \, 
\max(a_n, b_n).
\end{equation}
\end{theorem}
In most interesting situations $a_n$ is less than $b_n$ which allows
slightly to reformulate \eqref{eq:moment-bound}. Note that the
classical Stirling approximation to $n!$ gives
$$
\frac{(2n)!}{n! \, 2^n} \sim 2^{n+1/2} 
\, e^{-n} \, n^n \le  2^{n+1} \, e^{-n} \, n^n 
$$
and, therefore, letting $\phi$ be the `golden ratio'
$\phi := (1 + \sqrt{5})/2$, the $2n$th moment obeys 
\begin{equation}
    \label{eq:moment-bound2}
    \E(\tr(H_0^{2n})) \le 2 \, e^{-n} \, \gamma^{2n} \, n^{n+1} 
\cdot |\tau N|^n \, |T|^{n+1},  
\quad \gamma^2 = \frac{2 \phi^2}{1-\tau}, 
  \end{equation}
  provided that $a_n$ obeys
\begin{equation}
  \label{eq:n-small}
  a_n \le  2^{n+1} \, e^{-n} \, n^n \, 
\left(\frac{\tau}{1-\tau}\right)^n \, N^n \, |T|^n.  
\end{equation}

\newcommand{\goto}{\rightarrow}

Theorem \ref{moment-bound} gives a precise estimate about the operator
norm of $H_0$. To see why this is true, assume that \eqref{eq:n-small}
holds; since $H_0$ is self-adjoint
$$
\|H_0\|^{2n} = \|H_0^{n}\|^2 \le \|H_0^n\|^2_F = \tr(H_0^{2n})
$$
and, therefore,
\[
\left(\E \|H_0\|\right)^{2n} \le \E \|H_0\|^{2n} \le (2 n) \,
\gamma^{2n} \, e^{-n} \, n^{n} \, |T|^{n+1} \, |\tau N|^n.
\]
Now selecting $n = \lceil \log |T| \rceil$ so that 
$$
 e^{-n} \, n^{n} \, |T|  \le  \lceil \log |T| \rceil^n
$$
gives 
\[
\E \|H_0\| \le \gamma \cdot \sqrt{\log(|T|)} \cdot \sqrt{|T|\, |\tau
  N|} \cdot (1 + o(1)), \quad \text{ as } |T| \goto \infty.
\]
Formalizing matters, we proved 
\begin{corollary}
  Suppose $|T| \le (\log |\tau N|)^{-1} |\tau N|$. Then for any
  $\epsilon > 0$, we have
  $$
  \P\left(\|H_0\| > (1 + \epsilon) \, \gamma \, \sqrt{\log|T|} \,
    \sqrt{|T|\, |\tau N|}\right) \, \goto 0 \quad \text{ as } |T|,
  |\tau N| \goto \infty.
$$
\end{corollary}
\begin{proof} The Markov inequality above bounds the probability by
  $(1+\epsilon)^{-2n}$ which goes to zero as $n = \lceil \log |T|
  \rceil$ goes to infinity.
\end{proof}
  
We now return to the study of the invertibility of $\iota^* \iota -
\frac{1}{|\Omega|} H_0$. Letting $\alpha$ be a positive number $0 <
\alpha < 1$, it follows from the Markov inequality that
\begin{eqnarray*}
  \P( \|H_0^n\|_F \ge \alpha^n \cdot |\tau N|^n) 
= \frac{\E \|H_0^{n}\|^2_F}{\alpha^{2n} \, |\tau N|^{2n}}.
\end{eqnarray*}
We then apply inequality \eqref{eq:moment-bound} (recall
$\|H_0^{n}\|^2_F = \tr(H_0^{2n})$) and obtain
\begin{equation}
  \label{eq:large-dev}
   \P( \|H_0^n\|_F \ge \alpha^n \cdot |\tau N|^n) 
\le (2n) \,  e^{-n} \, \left(\frac{n \, \gamma^2}{\alpha^2}\right)^n \, 
\left(\frac{|T|}{|\tau N|}\right)^n\ \, |T|. 
\end{equation}
We remark that the last inequality holds for any sample size $|T|$
(proviso the condition \eqref{eq:n-small}) and we
now specialize \eqref{eq:large-dev} to selected values of $|T|$.

Suppose that $|T|$ obeys
\begin{equation}
  \label{eq:goodE}
|T| \le  \frac{\alpha_M^2}{\gamma^2} \,
\frac{|\tau N|}{n} < |T| + 1, \quad \text{for some } \alpha_M \le \alpha.  
\end{equation}
Then 
$$
\P( \|H_0^n\|_F \ge \alpha^n \cdot |\tau N|^n) \le 2 (\alpha^2/\gamma^2) \, 
e^{-n} \, |\tau N|.
$$
We then have the following result. 
\begin{theorem}
\label{super-useful}
Assume that $\tau \le .44$, say, and suppose that $T$ obeys
\eqref{eq:goodE}. Then \eqref{eq:n-small} holds for any $n \ge 4$, and
therefore
  \begin{equation}
    \label{eq:useful}
\P( \|H_0^n\|_F \ge \alpha^n \cdot |\tau N|^n) \le 2 (\alpha/\gamma)^2 \, 
e^{-n} \, |\tau N|.
  \end{equation}
\end{theorem}
The only thing to establish is that $T$ obeys \eqref{eq:n-small}. This is 
merely technical and the proof is in the Appendix. 

With the notations of the previous section and especially
\eqref{eq:tail-binomial2}, observe now that
$$
\P( \|H_0\| \ge \alpha \cdot |\Omega|) \le \P(\|H_0\| \ge \alpha
\, (1 - \epsilon_M) |\tau N|) + \P(|\Omega| < (1 - \epsilon_M) |\tau N|),
$$
where we recall that $B_M := \{|\Omega| < (1 - \epsilon_M) |\tau
N|\}$ has probability less than $N^{-M}$. Suppose $T$ obeys
\eqref{eq:goodE} with $\alpha_M := \alpha (1 - \epsilon_M)$ instead of
$\alpha$,
$$
\P(\|H_0\| \ge \alpha \, (1 - \epsilon_M) 
\cdot |\tau N|) \le 2 (\alpha/\gamma)^2 \, e^{-n} \, |\tau N|.
$$
\begin{corollary}
Take $n = (M+1)\log N$. We see from the Neumann series that the
operator $\iota^* \iota - \frac{1}{|\Omega|} \iota^* H$ is invertible
with probability at least $1 - (1 + 2/\gamma^2) N^{-M}$ since $\iota^*
\iota$ is the identity on vectors supported on $T$.
\end{corollary}
We have thus established the invertibility of $\iota^* \iota -
\frac{1}{|\Omega|} \iota^* H$ with high probability, and thus $P$ is
well defined with high probability.  It remains to show that $\sup_{t
  \notin T} |P(t)| < 1$ with high probability.

\subsection{Magnitude of the polynomial on the complement of $T$}
  
We first develop an expression for $P(t)$ by making use of the
algebraic identity
$$
(1-M)^{-1} = (1-M^n)^{-1} (1 + M + \ldots + M^{n-1}). 
$$
Indeed, we can write
$$
(\iota^* \iota - \frac{1}{|\Omega|^n} (\iota^*
    H)^n)^{-1} = \iota^* \iota + R$$
so that the inverse
is given by the truncated Neumann series
\begin{equation}\label{resolvent}
 (\iota^* \iota - \frac{1}{|\Omega|} \iota^* H)^{-1} 
= (\iota^* \iota + R)\sum_{m=0}^{n-1} \frac{1}{|\Omega|^m} (\iota^* H)^m. 
\end{equation}

The point is that the remainder term $R$ is quite small in the
Frobenius norm: suppose that $\|\iota^* H\|_F \le \alpha \cdot
|\Omega|$, then
$$
\|R\|_F \le \frac{\alpha^{n}}{1 - \alpha^n}.  
$$
In particular, the matrix coefficients of $R$ are all individually
less than $\alpha^n/(1 - \alpha^n)$. Introduce the $\ell_\infty$-norm
of a matrix as $\|M\|_\infty = \sup_{\|x\|_\infty \le 1} \|M
x\|_\infty$ which is also given by
$$
\|M\|_{\infty} = \sup_i \sum_j |M(i,j)|.
$$
Now, it follows from the Cauchy-Schwarz inequality that 
$$
\|M\|^2_{\infty} \le \sup_i \, \# M(\text{col}) \, \sum_j |M(i,j)|^2 \le
\# M(\text{col}) \cdot \|M\|^2_F,
$$
where $\# \, M(\text{col})$ is of course the number of columns of
$M$. This observation gives the crude estimate
\begin{equation}\label{r-small}
\| R\|_{\infty} \leq |T|^{1/2} \cdot \frac{\alpha^n}{1- \alpha^n}.
\end{equation}
As we shall soon see, the bound \eqref{r-small} allows us to
effectively neglect the $R$ term in this formula; the only remaining
difficulty will be to establish good bounds on the truncated Neumann
series $\frac{1}{|\Omega|} H \sum_{m = 0}^{n-1} \frac{1}{|\Omega|^m}
(\iota^* H)^m$.

\subsection{Estimating the truncated Neumann series}

From \eqref{P-def} we observe that on the complement of $T$ 
$$
P = \frac{1}{|\Omega|}H (\iota^* \iota -
  \frac{1}{|\Omega|}\iota^* H)^{-1} \iota^* \sgn(f), 
  $$
  since the $\iota$ component in \eqref{P-def} vanishes outside of
  $T$.  Applying \eqref{resolvent}, we may rewrite $P$ as
$$
P(t) = P_0(t) + P_1(t), \qquad \forall t \in T^c, 
$$
where
$$
P_0 = S_n \sgn(f),\quad P_1 = \frac{1}{|\Omega|} HR \iota^* (I +
S_{n-1}) \sgn(f)
$$
and 
$$
S_n = \sum_{m = 1}^{n} |\Omega|^{-m} (H \iota^*)^m.
$$
Let $a_0, a_1 > 0$ be two numbers with $a_0 + a_1 = 1$. Then
\[
\P\left(\sup_{t \in T^c} |P(t)| > 1\right) \le \P( \|P_0\|_\infty > a_0)
+ \P(\|P_1\|_\infty > a_1), 
\]
and the idea is to bound each term individually. Put $Q_0 = S_{n-1}
\sgn(f)$ so that $P_1 = \frac{1}{|\Omega|} H R \iota^* (\sgn(f) +
Q_0)$. With these notations, observe that
$$
\|P_1\|_\infty \le \frac{1}{|\Omega|} \|HR\|_\infty (1 +
\|\iota^*Q_0\|_\infty).
$$
Hence, bounds on the magnitude of $P_1$ will follow from bounds on
$\|H R\|_\infty$ together with bounds on the magnitude of $\iota^*
Q_0$.  It will be of course sufficient to derive bounds on
$\|Q_0\|_\infty$ (since $\|\iota^* Q_0\|_\infty \le \|Q_0\|_{\infty}$)
which will follow from those on $P_0$ since $Q_0$ is nearly equal to
$P_0$ (they differ by only one very small term term). 

Fix $t \in T^c$ and write $P_0(t)$ as 
$$
P_0(t) = \sum_{m = 1}^{n} |\Omega|^{-m} X_m(t), \qquad X_m = (H
\iota^*)^m \, \sgn(f)
$$
The idea is to use moment estimates to control the size of each
term $X_m(t)$.
\begin{lemma}
\label{second-moment}
Set $n = k m$. Then $\E|X_m(t_0)|^{2k}$ obeys the same estimate as
that in Theorem \ref{moment-bound} (up to a multiplicative factor
$|T|^{-1}$), namely,
\begin{equation}
  \label{eq:second-moment}
  \E|X_m(t_0)|^{2k} \le \frac{1}{|T|}
\cdot   n \, \phi^{2n} \,  \max(a_n, b_n). 
\end{equation}
In particular, following \eqref{eq:moment-bound2}
\begin{equation}
    \label{eq:second-moment-bound}
    \E|X_m(t_0)|^{2k}  \le 2 \, e^{-n} \, \gamma^{2n} \, n^{n+1} 
\cdot |T|^n |\tau N|^{n}, 
  \end{equation}
  where $\gamma$ is as before.
\end{lemma}
The proof of these moment estimates mimics that of Theorem
\ref{moment-bound} and may be found in the Appendix.

\begin{lemma}
\label{teo:P0(t)}
Fix $a_0 = .91$.  Suppose that $|T|$ obeys \eqref{eq:goodE} and let
$B_M$ be the set where $|\Omega| < (1-\epsilon_M) \cdot |\tau N|$
with $\epsilon_M$ as in \eqref{eq:tail-binomial2}.  For each $t \in
\ZZ_N$, there is a set $A_t$ with the property
\[
\P(A_t) > 1 - \epsilon_n, \qquad \epsilon_n = 2 (1 - \epsilon_M)^{-2n}
\cdot n^2 \, e^{-n} \alpha^{2n} \cdot (0.42)^{-2n},
\]
and 
\[
|P_0(t)| < .91, \quad |Q_0(t)| < .91 \text{ on } A_t \cap B_M^c.
\]
As a consequence, 
\[
\P(\sup_t |P_0(t)| > a_0) \le N^{-M} + N \epsilon_n,
\]
and similarly for $Q_0$.
\end{lemma}
\begin{proof}
  We suppose that $n$ is of the form $n = 2^J - 1$ (this property is
  not crucial and only simply simplifies our exposition). For each $m$
  and $k$ such that $k m \ge n$, it follows from \eqref{eq:goodE} and
  \eqref{eq:second-moment-bound} together with some simple
  calculations that 
\begin{equation}
  \label{eq:simple}
  \E |X_m(t)|^{2k} \le 
2 n \, e^{-n} \alpha^{2n} \cdot |\tau N|^{2n}.
\end{equation}

Again $|\Omega| \approx |\tau N|$ and we will develop a bound on the
set $B_M^c$ where $|\Omega| \geq (1-\epsilon_M) |\tau N|$. On this set 
\[
|P_0(t)| \le \sum_{m = 1}^{n} Y_m, \qquad Y_m =
\frac{1}{(1 -\epsilon_M)^m \, |\tau N|^m } \, |X_m(t)|.
\]

Fix $\beta_j > 0$, $0 \le j < J$, such that $\sum_{j = 0}^{J-1}
2^j \, \beta_j \le a_0$.  Obviously, 
$$
\P(\sum_{m = 1}^n Y_m > a_0) \le \sum_{j = 0}^{J-1} \,\, \sum_{m =
  2^j}^{2^{j+1} - 1} \P(Y_m > \beta_j) \le \sum_{j = 0}^{J-1} \,\,
\sum_{m = 2^j}^{2^{j+1} - 1} \beta_j^{-2K_j} \, \E|Y_m|^{2K_j}. 
$$
where $K_j = 2^{J-j}$. Observe that for each $m$ with $2^j \le m <
2^{j+1}$, $K_j m$ obeys $n \le K_j m < 2n$ and, therefore,
\eqref{eq:simple} gives
$$
\E|Y_m|^{2K_j} \le (1-\epsilon_M)^{-2n} \cdot (2 n \, e^{-n}
\alpha^{2n}).
$$
For example, taking $\beta_j^{-K_j}$ to be constant for all $j$,
i.e.~equal to $\beta_0^{-n}$, gives
$$
\P(\sum_{m=1}^n Y_m > a_0) \le 2 (1-\epsilon_M)^{-2n} \cdot n^2 \, e^{-n}
  \alpha^{2n} \cdot \beta_0^{-2n},  
  $$
  with $\sum_{j = 0}^{J-1} 2^j \beta_j \le a_0$.  Numerical
  calculations show that for $\beta_0 = .42$, $\sum_j 2^j \beta_j \le
  .91$ which gives
\begin{equation}
  \label{eq:best-beta}
  \P(\sum_{m=1}^n Y_m > .91) \le  2 (1-\epsilon_M)^{-2n} \cdot n^2 \, e^{-n}
  \alpha^{2n} \cdot (0.42)^{-2n}. 
\end{equation}
The claim for $Q_0$ is, of course, identical and the lemma follows.
\end{proof} 
\begin{lemma}
\label{teo:P1(t)}
Fix $a_1 = .09$. Suppose that the pair $(\alpha, N)$ obeys $|\tau N|^{3/2}
\frac{\alpha^n}{1 - \alpha^n} \le a_1/2$. Then
$$
\|P_1\|_\infty \le a_1
$$
on the event $A \cap \{\|\iota^* H\|_F \le \alpha |\Omega|\}$, for
some $A$ obeying $\P(A) \ge 1 - O(N^{-M})$.
\end{lemma}
\begin{proof}
  As we observed before, (1) $\|P_1\|_\infty \le \|H\|_\infty
  \|R\|_\infty (1 + \|Q_0\|_\infty)$, and (2) $Q_0$ obeys the bound
  stated in Lemma \ref{teo:P0(t)}. Consider then the event
  $\{\|Q_0\|_\infty \le 1\}$. On this event, $\|P_1\| \le a_1$ if
  $\frac{1}{|\Omega|} \|H\| \|R\|_\infty \le a_1/2$.  The matrix $H$ obeys
  $\frac{1}{|\Omega|} \|H\|_{\infty} \le |T|$ since $H$ has $|T|$
  columns and each matrix element is bounded by $|\Omega|$ (note that
  far better bounds are possible). It then follows from
  \eqref{r-small} that
  $$
  \|H\|_{\infty} \cdot \|R\|_{\infty} \le |T|^{3/2} \cdot
  \frac{\alpha^n}{1 - \alpha^n},
$$
with probability at least $1 - O(N^{-M})$. We then simply need to choose 
$\alpha$ and $n$ such that the right hand-side is less than $a_1/2$. 
\end{proof}

\subsection{Proof of Theorem \ref{main}}
\label{sec:proof}

It is now clear that we have assembled all the intermediate results to
prove our theorem. Indeed, we proved the invertibility of $i^*i -
\frac{1}{|\Omega|} \iota^* H$ with probability $O(N^{-M})$ and $|P(t)|
< 1$ for all $t \in T^c$ (again with high probability), provided that
$\alpha$ and $n$ be selected appropriately as we now explain.

Fix $M > 0$. We choose $\alpha = .42$ and $n$ to be the nearest
integer to $(M+1) \log N$.
\begin{enumerate}
\item From the discussion following Theorem \ref{super-useful}, it
  follows that $i^*i - |\Omega|^{-1} \iota^* H$ is invertible with
  probability $O(N^{-M})$.
\item With this special choice, $\epsilon_n = 2[(M+1) \log N]^2 \cdot
  N^{-(M+1)}$ and, therefore, Lemma \ref{teo:P0(t)} implies that both
  $P_0$ and $Q_0$ are bounded by .91 outside of $T^c$ with probability
  at least $1 - [1 + 2((M+1) \log N)^2] \cdot N^{-M}$.
\item And finally, to prove that $|P_1(t)| < .09$ outside $T^c$, Lemma
  \ref{teo:P0(t)} assures that it is sufficient to have $N^{3/2}
  \alpha^n/(1-\alpha^n) \le .045$. Because $\log(.42) \approx -.87$
  and $\log(.045) \approx -3.10$, this condition is approximately
  equivalent to
\[
(1.5 - .87 (M+1)) \log N \le -3.10.
\]
Take $M \ge 2$, for example; then the above inequality is satisfied as
soon as $N \ge 17$.
\end{enumerate}
To conclude, we proved that if $T$ obeys 
\[
|T| \le \alpha(M) \cdot \frac{|\tau N|}{\log N}, \qquad \alpha(M) =
\frac{.42^2}{\gamma^2 (M+1)}(1 + o(1))
\]
then the reconstruction with probability exceeding $1 - O([(M+1) \log
N)^2] \cdot N^{-M})$. In other words, we may take $\alpha(M)$ in
Theorem \ref{main} to be of the form
\begin{equation}
  \label{eq:final-alpha}
\alpha(M) = \frac{1}{29.6(M+1)} (1 + o(1)). 
\end{equation}

\section{Moments of Random Matrices}
\label{sec:moments}

\subsection{A First Formula for the Expected Value 
  of the Trace of $(H_0)^{2n}$}

Recall that $H_0(t,t')$, $t, t' \in T$, is the $|T| \times |T|$
matrix whose entries are defined by
\begin{equation}
  \label{eq:T}
H_0(t,t') =  \begin{cases} 0 & t = t',\\ 
 c(t - t') & t \neq t',    \end{cases}  
\quad c(u) = \sum_{\omega \in \Omega} e^{i \omega u}. 
\end{equation}
A diagonal element of the $2n$th power of $H_0$ may be expressed as
$$
H_0^{2n}(t_1,t_1) = \sum_{t_2, \ldots, t_{2n}: \, t_j \neq t_{j+1}}
\, c(t_1 - t_2) \ldots c(t_{2n} - t_1), 
$$
where we adopt the convention that $t_{2n+1} = t_1$ whenever
convenient and, therefore,
$$
\E(\tr(H_0^{2n})) = \sum_{t_1, \ldots, t_{2n}: \, t_j \neq t_{j+1}}
\,\,\, \E\left[ \sum_{\omega_1, \ldots, \omega_{2n} \in \Omega} \,\,\,
  e^{i \sum_{j=1}^{2n} \omega_j (t_j - t_{j+1})} \right].
$$
Using \eqref{omega-def} and linearity of expectation, we can write
this as
$$
\sum_{t_1, \ldots, t_{2n}: \, t_j \neq t_{j+1}} \,\,\,
\sum_{0 \le \omega_1, \ldots, \omega_{2n} \le N-1} \,\,\,
e^{i \sum_{j=1}^{2n} \omega_j (t_j - t_{j+1})} \, \E \left[\prod_{j=1}^{2n}
I_{\{\omega_j \in \Omega\}}\right].
$$

The idea is to use the independence of the $I_{\{\omega_j \in \Omega\}}$'s
to simplify this expression substantially; however, one has to be
careful with the fact that some of the $\omega_j$'s may be the same,
at which point one loses independence of those indicator variables.
These difficulties require a certain amount of notation.  We let
$\ZZ_N = \{0, 1, \ldots, N - 1\}$ be the set of all frequencies as
before, and let $A$ be the finite set $A := \{1, \ldots, 2n\}$. For
all $\bomega := (\omega_1, \ldots, \omega_{2n})$, we define the
equivalence relation $\sim_{\bomega}$ on $A$ by saying that $j
\sim_\bomega j'$ if and only if $\omega_j = \omega_{j'}$. We let ${\cal
  P}(A)$ be the set of all equivalence relations on $A$.  Note that
there is a partial ordering on the equivalence relations as one can
say that $\sim_1 \leq \sim_2$ if $\sim_1$ is coarser than $\sim_2$,
i.e.  $a \sim_2 b$ implies $a \sim_1 b$ for all
$a,b \in A$.  Thus, the coarsest element in ${\cal P}(A)$ is
the trivial equivalence relation in which all elements of $A$ are
equivalent (just one equivalence class), while the finest element is
the equality relation $=$, i.e. each element of $A$ belongs to
a distinct class ($|A|$ equivalence classes).

For each equivalence relation $\sim$ in ${\cal P}$, we can then define
the sets $\Omega(\sim) \subset \ZZ_N^{2n}$ by
$$
\Omega(\sim) := \{\bomega \in \ZZ_N^{2n}: \sim_\bomega = \sim \}$$
and the
sets $\Omega_{\leq}(\sim) \subset \ZZ_N^{2n}$ by
$$
\Omega_{\leq}(\sim) := \bigcup_{\sim' \in {\cal P}: \sim' \leq
  \sim} \Omega(\sim') = \{ \bomega \in \ZZ_N^{2n}: \sim_\bomega \leq \sim \}.
$$
Thus the sets $\{ \Omega(\sim): \sim \in {\cal P}\}$ form a
partition of $\ZZ_N^{2n}$.  The sets $\Omega_{\leq}(\sim)$ can also be
defined as
$$
\Omega_{\leq}(\sim) := \{ \bomega \in \ZZ_N^{2n}: \omega_a =
\omega_b \hbox{ whenever } a \sim b \}.
$$
For comparison, the sets $\Omega(\sim)$ can be defined as
$$
\Omega(\sim) := \{ \bomega \in \ZZ_N^{2n}: \omega_a = \omega_b
\hbox{ whenever } a \sim b, \hbox{ and } \omega_a \neq
\omega_b \hbox{ whenever } a \not \sim b\}.
$$
We give an example: suppose $n = 2$ and fix $\sim$ such that $1
\sim 4$ and $2 \sim 3$ (exactly 2 equivalence classes); then
$\Omega(\sim) := \{\bomega \in \ZZ_N^4 : \omega_1 = \omega_4, \, \omega_2 =
\omega_3, \text{ and } \omega_1 \neq \omega_2\}$ while $\Omega_\leq(\sim)
:= \{\bomega \in \ZZ_N^4 : \omega_1 = \omega_4, \, \omega_2 = \omega_3\}$.

Now, let us return to the computation of the expected value.  Because
the random variables $I_k$ \eqref{eq:Ik} are independent and have all
the same distribution, the quantity $\E[\prod_{j=1}^{2n}
I_{\omega_j}]$ depends only on the equivalence relation
$\sim_{\bomega}$ and not on the value of $\bomega$ itself.  Indeed, we
have
$$
\E(\prod_{j=1}^{2n} I_{\omega_j}) = \tau^{|A/\sim|},$$
where
$A/\sim$ denotes the equivalence classes of $\sim$.  Thus we can
rewrite the preceding expression as
\begin{equation}
  \label{eq:formula1}
\E(\tr(H_0^{2n})) = 
\sum_{t_1, \ldots, t_{2n}: \, t_j \neq t_{j+1}}
\,\,\, \sum_{\sim \in {\cal P}(A)} \,\, \tau^{|A/\sim|} 
\sum_{\bomega \in \Omega(\sim)} 
\,\,\, e^{i \sum_{j=1}^{2n} \omega_j (t_j - t_{j+1})}
\end{equation}
where $\sim$ ranges over all equivalence relations.

We would like to pause here and consider \eqref{eq:formula1}. Take $n
= 1$, for example. There are only two equivalent classes on $\{1,2\}$
and, therefore, the right hand-side is equal to
\begin{equation*}
\sum_{t_1, t_2: \, t_1 \neq t_2}
\,\,\, \left[\tau \sum_{(\omega_1, \omega_2) \in \ZZ_N^2 : \omega_1 = \omega_2} 
e^{ i\omega_1 (t_1 - t_1)} + \tau^2
\sum_{(\omega_1, \omega_2) \in \ZZ_N^2 : \omega_1 \neq \omega_2} 
e^{ i \omega_1 (t_1 - t_2) + i \omega_2(t_2 - t_1)}\right]. 
\end{equation*}
Our goal is to rewrite the expression inside the brackets so that the
exclusion $\omega_1 \neq \omega_2$ does not appear any longer, i.e.
we would like to rewrite the sum over $\bomega \in \ZZ_N^2 : \omega_1 \neq
\omega_2$ in terms of sums over $\bomega \in \ZZ_N^2 : \omega_1 =
\omega_2$, and over $\bomega \in \ZZ_N^2$. In this special case, this is
quite easy as 
$$
\sum_{\bomega \in \ZZ_N^2 : \omega_1 \neq \omega_2} = \sum_{\bomega \in
  \ZZ_N^2} - \sum_{\bomega \in \ZZ_N^2 : \omega_1 = \omega_2}
$$
The motivation is quite clear. Removing the exclusion allows to
rewrite sums as product, e.g.
$$
\sum_{\bomega \in \ZZ_N^2} = \sum_{\omega_1} e^{ i \omega_1 (t_1 - t_2)}
\cdot \sum_{\omega_2} e^{ i \omega_2 (t_2 - t_1)}; 
$$
and each factor is equal to either $N$ or $0$ depending on whether
$t_1 = t_2$ or not.

The next section generalizes these ideas and develop an identity,
which allows us to rewrite sums over $\Omega(\sim)$ in terms of sums
over $\Omega_{\leq}(\sim)$.

\subsection{Inclusion-Exclusion formulae}

\begin{lemma}[Inclusion-Exclusion principle for equivalence classes]  
  Let $A$ and $G$ be non-empty finite sets.  For any equivalence class
  $\sim \in {\cal P}(A)$ on $\bomega\in G^{|A|}$, we have
\begin{equation}\label{inclusion}
 \sum_{\bomega \in \Omega(\sim)} f(\bomega) = 
\sum_{\sim_1 \in {\cal P}: \sim_1 \leq \sim}
(-1)^{|A/\sim| - |A/\sim_1|} 
\left(\prod_{A' \in A/\sim_1} (|A'/\sim|-1)!\right) 
\sum_{\bomega \in \Omega_{\leq}(\sim_1)} f(\bomega).
\end{equation}
\end{lemma}
Thus, for instance, if $A = \{1,2,3\}$ and $\sim$ is the equality
relation, i.e. $j \sim k$ if and only if $j = k$, this identity is
saying that
\begin{multline*}
  \sum_{\omega_1,\omega_2,\omega_3 \in G: \omega_1,\omega_2,\omega_3
  \text{ distinct}} = \sum_{\omega_1,\omega_2,\omega_3 \in G}\\
 - \sum_{\omega_1,\omega_2,\omega_3: \omega_1 = \omega_2} -
\sum_{\omega_1,\omega_2,\omega_3 \in G: \omega_2 = \omega_3} -
\sum_{\omega_1,\omega_2,\omega_3 \in G: \omega_3 = \omega_1} +
2\sum_{\omega_1,\omega_2,\omega_3 \in G: \omega_1 = \omega_2 =
  \omega_3}
\end{multline*}
where we have omitted the summands
$f(\omega_1,\omega_2,\omega_3)$ for brevity.

\begin{proof} 
  By passing from $A$ to the quotient space $A/\sim$ if necessary we
  may assume that $\sim$ is the equality relation $=$.  Now relabeling
  $A$ as $\{1,\ldots,n\}$, $\sim_1$ as $\sim$, and $A'$ as $A$, it
  suffices to show that
\begin{multline}\label{G-induct}
 \sum_{\bomega \in G^n: \omega_1, \ldots, \omega_n 
\hbox{ distinct}} f(\bomega) = \\
\sum_{\sim \in {\cal P}(\{1,\ldots,n\})}
(-1)^{n - |\{1,\ldots,n\}/\sim|} 
\left[\prod_{A \in \{1,\ldots,n\}/\sim} (|A|-1)!\right] 
\sum_{\bomega \in \Omega_{\leq}(\sim)} f(\bomega).
\end{multline}
We prove this by induction on $n$.  When $n=1$ both sides are equal to
$\sum_{\bomega \in G} f(\bomega)$.  Now suppose inductively that $n > 1$
and the claim has already been proven for $n-1$.  We observe that the
left-hand side of \eqref{G-induct} can be rewritten as
$$
\sum_{\bomega' \in G^{n-1}: \omega_1, \ldots, \omega_{n-1}
  \hbox{ distinct}} \left(\sum_{\omega_n \in G} f(\bomega',\omega_n) -
\sum_{j=1}^{n-1} f(\bomega',\omega_j)\right), 
$$
where $\bomega' := (\omega_1,\ldots,\omega_{n-1})$.  Applying the
inductive hypothesis, this can be written as
\begin{multline}\label{G-LHS}
\sum_{\sim' \in {\cal P}(\{1,\ldots,n-1\})}
(-1)^{n-1 - |\{1,\ldots,n-1\}/\sim'|} \, 
\prod_{A' \in \{1,\ldots,n-1\}/\sim} (|A'|-1)! \\
\sum_{\bomega' \in \Omega_{\leq}(\sim')} 
\left(\sum_{\omega_n \in G} f(\bomega',\omega_n) - 
\sum_{1 \le j \le n} f(\bomega',\omega_j)\right).
\end{multline}
Now we work on the right-hand side of \eqref{G-induct}.  If $\sim$ is
an equivalence class on $\{1,\ldots,n\}$, let $\sim'$ be the
restriction of $\sim$ to $\{1,\ldots,n-1\}$.  Observe that $\sim$ can
be formed from $\sim'$ either by adjoining the singleton set $\{n\}$
as a new equivalence class (in which case we write $\sim = \{ \sim',
\{n\} \}$, or by choosing a $j \in \{1,\ldots,n-1\}$ and declaring $n$
to be equivalent to $j$ (in which case we write $\sim = \{ \sim',
\{n\} \}/(j=n)$).  Note that the latter construction can recover the
same equivalence class $\sim$ in multiple ways if the equivalence
class $[j]_{\sim'}$ of $j$ in $\sim'$ has size larger than 1, however
we can resolve this by weighting each $j$ by
$\frac{1}{|[j]_{\sim'}|}$.  Thus we have the identity
\begin{multline*}
\sum_{\sim \in {\cal P}(\{1,\ldots,n\})} F(\sim) = \sum_{\sim' \in
  {\cal P}(\{1,\ldots,n-1\})} F(\{\sim', \{n\}\}) \\
+ \sum_{\sim' \in
  {\cal P}(\{1,\ldots,n-1\})} \sum_{j=1}^{n-1} \frac{1}{|[j]_{\sim'}|}
F(\{\sim', \{n\}\}/(j=n))
\end{multline*}
for any complex-valued function $F$ on
${\cal P}(\{1,\ldots,n\})$.  Applying this to the right-hand side of
\eqref{G-induct}, we see that we may rewrite this expression as the
sum of
$$\sum_{\sim' \in {\cal P}(\{1,\ldots,n-1\})} (-1)^{n -
  (|\{1,\ldots,n-1\}/\sim'|+1)} \left[\prod_{A \in \{1,\ldots,n-1\}/\sim'}
(|A|-1)! \right]
\sum_{\bomega' \in \Omega_{\leq}(\sim')} f(\bomega',\omega_n)$$
and
$$
\sum_{\sim' \in {\cal P}(\{1,\ldots,n-1\})}  
(-1)^{n - |\{1,\ldots,n-1\}/\sim'|}  \,\,\, 
\sum_{j=1}^{n-1} T(j)  \sum_{\bomega'
\in \Omega_{\leq}(\sim')} f(\bomega',\omega_j),
$$
where we adopt the convention $\bomega' = (\omega_1,\ldots,\omega_{n-1})$.
But observe that
$$
T(j): = \frac{1}{|[j]_{\sim'}|} \prod_{A \in
  \{1,\ldots,n\}/(\{\sim',\{n\}\}/(j=n))} (|A|-1)!  = \prod_{A' \in
  \{1,\ldots,n-1\}/\sim'} (|A'|-1)!$$
and thus the right-hand side of
\eqref{G-induct} matches \eqref{G-LHS} as desired.
\end{proof}

\subsection{Stirling Numbers}

As emphasized earlier, our goal is to use our inclusion-exclusion
formula to rewrite the sum \eqref{eq:formula1} as a sum over
$\Omega_{\leq}(\sim)$. In order to do this, it is best to introduce
another element of combinatorics, which will prove to be very useful.

For any $n, k \geq 0$, we define the \emph{Stirling number of the
  second kind} $S(n,k)$ to be the number of equivalence relations on a
set of $n$ elements which have exactly $k$ equivalence classes, thus
$$
S(n,k) := \#\, \{ \sim \in {\cal P}(A): |A/\sim| = k \}. 
$$
Thus for instance $S(0,0) = S(1,1) = S(2,1) = S(2,2) = 1$, $S(3,2)
= 3$, and so forth.  We observe the basic recurrence
\begin{equation}\label{recurrence}
 S(n+1,k) = S(n,k-1) + k S(n,k) 
\text{ for all } k, n \geq 0. 
\end{equation}
This simply reflects the fact that if $a$ is an element of $A$
and $\sim$ is an equivalence relation on $A$ with $k$ equivalence
classes, then either $a$ is not equivalent to any other element
of $A$ (in which case $\sim$ has $k-1$ equivalence classes on $A
\backslash \{a\}$), or $a$ is equivalent to one of the $k$
equivalence classes of $S \backslash \{a\}$.

We now need an identity for the Stirling numbers\footnote{We found
  this identity by modifying a standard generating function identity
  for the Stirling numbers which involved the polylogarithm.  It can
  also be obtained from the formula $S(n,k) = \frac{1}{k!}
  \sum_{i=0}^{k-1} (-1)^i {k \choose i} (k-i)^n$, which can be
  verified inductively from \eqref{recurrence}.}.
\begin{lemma}\label{poly}
For any $n \geq 1$ and $0 \leq \tau < 1/2$, we have the identity
\begin{equation}\label{tau-ident}
\sum_{k=1}^n (k-1)! S(n,k) (-1)^{n-k} \tau^k = \sum_{k=1}^\infty (-1)^{n-k} \frac{\tau^k k^{n-1}}{(1-\tau)^k}.
\end{equation}
Note that the condition $0 \leq \tau < 1/2$ ensures that the right-hand side is convergent.
\end{lemma}

\begin{proof} 
  We prove this by induction on $n$.  When $n=1$ the left-hand side is
  equal to $\tau$, and the right-hand side is equal to
  $$
  \sum_{k=1}^\infty (-1)^{k+1} \frac{\tau^k}{(1-\tau)^k} = -
  \sum_{k=0}^\infty \left(\frac{\tau}{\tau-1}\right)^k + 1 = \frac{-1}{1 -
    \frac{\tau}{\tau-1}} + 1 = \tau$$
  as desired.  Now suppose
  inductively that $n \geq 1$ and the claim has already been proven
  for $n$.  Applying the operator $(\tau^2 - \tau) \frac{d}{d\tau}$ to
  both sides (which can be justified by the hypothesis $0 \leq \tau <
  1/2$) we obtain (after some computation)
$$ \sum_{k=1}^{n+1} (k-1)! (S(n,k-1) + k S(n,k)) (-1)^{n+1-k} \tau^k  = 
\sum_{k=0}^\infty (-1)^{n+1-k} \frac{\tau^k k^{n}}{(1-\tau)^k},$$
and the claim follows from \eqref{recurrence}.
\end{proof}

We shall refer to the quantity in \eqref{tau-ident} as $F_n(\tau)$,
thus
\begin{equation}\label{tau-ident-n}
F_n(\tau) = \sum_{k=1}^n (k-1)! S(n,k) (-1)^{n-k} \tau^k = 
\sum_{k=1}^\infty (-1)^{n+k} \frac{\tau^k k^{n-1}}{(1-\tau)^k}.
\end{equation}
Thus we have 
$$
F_1(\tau) = \tau, \quad F_2(\tau) = -\tau + \tau^2, \quad 
F_3(\tau) = \tau - 3\tau^2 + 2 \tau^3,
$$
and so forth.  When $\tau$ is small we have the approximation
$F_n(\tau) \approx (-1)^{n+1} \tau$, which is worth keeping in mind.
Some more rigorous bounds in this spirit are as follows.
\begin{lemma}\label{fn-bound}  
  Let $n \geq 1$ and $0 \leq \tau < 1/2$.  If $\frac{\tau}{1-\tau}
  \leq e^{1-n}$, then we have $|F_n(\tau)| \leq \frac{\tau}{1-\tau}$.
  If instead $\frac{\tau}{1-\tau} > e^{1-n}$, then
  $$
  |F_n(\tau)| \leq \exp( (n-1) (\log(n-1) - \log \log
  \frac{1-\tau}{\tau} - 1) ).
$$
\end{lemma}

\begin{proof} 
  Elementary calculus shows that for $x > 0$, the function $g(x) =
  \frac{\tau^x x^{n-1}}{(1-\tau)^x}$ is increasing for $x < x_*$ and
  decreasing for $x > x^*$, where $x_* := (n-1)/\log
  \frac{1-\tau}{\tau}$.  If $\frac{\tau}{1-\tau} \leq e^{1-n}$, then
  $x_* \leq 1$, and so the alternating series $F_n(\tau) =
  \sum_{k=1}^\infty (-1)^{n+k} g(k)$ has magnitude at most $g(1) =
  \frac{\tau}{1-\tau}$.  Otherwise the series has magnitude at most
$$ g(x_*) = \exp( (n-1) (\log(n-1) - \log \log \frac{1-\tau}{\tau} - 1) )$$
and the claim follows.
\end{proof}

Roughly speaking, this means that $F_n(\tau)$ behaves like $\tau$ for
$n = O(\log [1/\tau])$ and behaves like $(n / \log [1/\tau])^n$ for $n
\gg \log [1/\tau]$.

\subsection{A Second Formula for the Expected Value 
    of the Trace of $H_0^{2n}$}
  
Let us return to \eqref{eq:formula1}. The inner sum of
  \eqref{eq:formula1} can be rewritten as
$$
  \sum_{\sim \in {\cal P}(A)} \tau^{|A/\sim|} \sum_{\bomega \in
    \Omega(\sim)} f(\bomega)
  $$
  with $f(\bomega) := e^{i \sum_{1 \le j \le 2n} \omega_j(t_j -
    t_{j+1})}$. We prove the following useful identity:
\begin{lemma}
\label{teo:combinatorics}
\begin{equation}\label{combinatorics}
\sum_{\sim\in{\cal P}(A)} \tau^{|A/\sim|}
\sum_{\bomega\in\Omega(\sim)} f(\bomega)
=  \sum_{\sim_1 \in {\cal P}(A)} \left[ 
\sum_{\bomega \in \Omega_{\leq}(\sim_1)} f(\bomega) \right]
 \prod_{A' \in A/\sim_1} F_{|A'|}(\tau).
\end{equation}
\end{lemma}
\begin{proof}
  Applying \eqref{inclusion} and rearranging, we may rewrite this as
$$
\sum_{\sim_1 \in {\cal P}(A)} T(\sim_1) \, \sum_{\bomega \in
    \Omega_{\leq}(\sim_1)} f(\bomega), 
$$
where
$$
T(\sim_1) = \sum_{\sim \in {\cal P}(A): \sim \geq \sim_1}
\tau^{|A/\sim|} (-1)^{|A/\sim| - |A/\sim_1|} \prod_{A' \in A/\sim_1}
(|A'/\sim|-1)!.
$$
Splitting $A$ into equivalence classes $A'$ of $A/\sim_1$, we
observe that
$$
T(\sim_1) = \prod_{A' \in A/\sim_1} \sum_{\sim' \in {\cal P}(A')}
\tau^{|A'/\sim'|} (-1)^{|A'/\sim'| - |A'|} (|A'/\sim'|-1)!;
$$
splitting $\sim'$ based on the number of equivalence classes 
$|A'/\sim'|$, we can write this as
$$
\prod_{A' \in A/\sim_1} \sum_{k=1}^{|A'|} S(|A'|,k) \tau^k
(-1)^{|A'|-k} (k-1)! = \prod_{A' \in A/\sim_1} F_{|A'|}(\tau)
$$
by \eqref{tau-ident-n}.  Gathering all this together, we have
proven the identity \eqref{combinatorics}.
\end{proof}

We specialize \eqref{combinatorics} to the function $f(\bomega) :=
\exp(i \sum_{1 \le j \le 2n} \omega_j (t_j - t_{j+1}))$ and obtain 
\begin{equation}
    \label{eq:expected1}
     \E[\tr(H_0^{2n})] =   \sum_{\sim \in {\cal P}(A)} \,\, 
\sum_{t_1, \ldots, t_{2n} \in T: \, t_j \neq t_{j+1}} \,\,
\sum_{\bomega \in \Omega_{\leq}(\sim)} 
e^{i \sum_{j = 1}^{2n} \omega_j(t_j - t_{j+1})}\,
\prod_{A' \in A/\sim}  F_{|A'|}(\tau).
  \end{equation}
We now compute
$$
I(\sim) = \sum_{\bomega \in \Omega_{\leq}(\sim)} e^{i \sum_{1 \le j
    \le 2n} \omega_j (t_j - t_{j+1})}.
$$
For every equivalence class $A' \in A/\sim$, let $t_{A'}$ denote
the expression $t_{A'} := \sum_{a \in A'} (t_a - t_{a+1})$, and let
$\omega_{A'}$ denote the expression $\omega_{A'} := \omega_a$ for
any $a \in A'$ (these are all equal since $\bomega \in
\Omega_{\leq}(\sim)$).  Then
$$
I(\sim) = \sum_{(\omega_{A'})_{A' \in A/\sim} \in \ZZ_N^{|A/\sim|}}
e^{\sum_{A' \in A/\sim} i \omega_{A'} t_{A'}} = \prod_{A' \in A/\sim}
\,\,\, \sum_{\omega_{A'} \in \ZZ_N} e^{i \omega_{A'} t_{A'}}.
$$
We now see the importance of \eqref{eq:expected1} as the inner sum
equals $|\ZZ_N| = N$ when $t_{A'} = 0$ and vanishes otherwise.  Hence,
we proved the following:
\begin{lemma}
\label{expected1}
For every equivalence class $A' \in A/\sim$, let 
\mbox{$t_{A'} := \sum_{a \in A'} (t_a - t_{a+1})$}. Then
\begin{equation}
  \label{eq:expected3}
\E[\tr(H_0^{2n})]   =   \sum_{\sim \in {\cal P}(A)} \,\, 
\sum_{t \in T^{2n}: \, t_j \neq t_{j+1} \text{ and }  t_{A'} = 0
\text{ for all } A'} N^{|A/\sim|} \prod_{A' \in A/\sim} F_{|A'|}(\tau).
\end{equation}
\end{lemma}
This formula will serve as a basis for all of our estimates.  In
particular, because of the constraint $t_j \neq t_{j+1}$, we see that
the summand vanishes if $A/\sim$ contains any singleton equivalence
classes. This means, in passing, that the only equivalence classes
which contribute to the sum obey $|A/\sim| \leq n$.

\subsection{A First Bound on $\E[\tr(H_0^{2n})]$}

Let $\sim$ be an equivalence which does not contain any singleton.
Then the following inequality holds
$$
\# \,\, \{t \in T^{2n}: \, t_{A'} = 0 \text{ for all
} A' \in {A/\sim} \} \le |T|^{2n - |A/\sim| + 1}. 
$$
To see why this is true, observe that as linear combinations of
$t_1, \ldots, t_{2n}$, the expressions $t_j - t_{j+1}$ are all
linearly independent of each other except for the constraint $\sum_{j
  = 1}^{2n} t_j - t_{j+1} = 0$.  Thus we have $|A/\sim|-1$ independent
constraints in the above sum, and so the number of $t$'s obeying the
constraints is bounded by $|T|^{2n-|A/\sim|+1}$.

All the equivalence classes in the sum \eqref{eq:expected3} are
without singletons as otherwise $t_{A'} \neq  0$. Thus, for $n, k
\geq 0$, we let $P(n,k)$ be the number of equivalence classes on a set
of $n$ elements which have exactly $k$ equivalence classes and {\em no
  singletons}
$$
P(n,k) := \# \, \{ \sim \in {\cal P}(A): \, |A/\sim| = k \text{ and
} |A'| \ge 2, \, \forall A' \in A/\sim\}.
$$
There is a simple recursion on these numbers, namely, 
\begin{equation}\label{recurrence2}
P(n,k) = P(n-1,k) + (n-1) P(n-2,k-1), 
\end{equation}
which is valid for all $n, k \ge 0$.  This simply reflects the fact
that if $\alpha$ is an element of $A$ and $\sim$ is an equivalence
relation on $A$ with $k$ equivalence classes, then either (1) $\alpha$
belongs to a class which has only one other element $\beta$ of $A$ (in
which case $\sim$ has $k-1$ equivalence classes and no singleton on $A
\backslash \{\alpha,\beta\}$), or $\alpha$ is equivalent to one of the
$k$ equivalence classes of $A \backslash \{\alpha\}$, each of which
having at least two elements.

With these notations, we established
\begin{equation}
  \label{eq:expected4}
 \E\tr(H^{2n}_0)   
\le   \sum_{k = 1}^{n}  N^{k} \, |T|^{2n - k + 1} \, P(2n,k) \, 
\sup_{\sim: |A/\sim| = k} \,\, \prod_{A' \in A/\sim} F_{|A'|}(\tau). 
\end{equation}
The following lemma provides an upper bound on those $P(n,k)$'s.
\begin{lemma}
 \label{Pnk}
  The numbers $P(n,k)$ obey 
  \begin{equation}
    \label{eq:Pnk}
    P(n,k) \le \lambda^{n} \, (n-1) \ldots (n - 2k + 1), 
\quad \forall \lambda \ge \phi := \frac{1 + \sqrt{5}}{2}.  
  \end{equation}
\end{lemma}
\begin{proof}
The proof operates by induction. The bound \eqref{eq:Pnk}
is obvious for $n = 1$. Suppose the claim is established for all pairs
$(m,k)$ with $m \le n$. We will show that this implies the property
for $m = n + 1$. Indeed,
\begin{eqnarray*}
  P(n+1,k) & = & P(n-1,k) + (n-1) P(n-2,k-1)\\
& \le & \lambda^{n-1} \, (n-2) \ldots (n - 2k + 2) +  \lambda^{n-2} 
\, (n-1) \ldots (n - 2k + 1)\\
& \le &  (\lambda^{n-1} + \lambda^{n-2}) \, (n-1) \ldots (n - 2k + 1). 
\end{eqnarray*}
The claim follows since for $\lambda \ge \frac{1 + \sqrt{5}}{2}$, we
have $\lambda^{n-1} + \lambda^{n-2} \le \lambda^{n}$. 
\end{proof}

This lemma gives us an idea of how large the $P(2n,k)$'s appearing in
the sum \eqref{eq:expected4} really are. To derive an upper bound on
the whole sum, we also need to understand the behavior of $ \prod_{A'
  \in A/\sim} F_{|A'|}(\tau)$. This is the subject of our next
section.

\subsection{Convex analysis}

We start with a useful and classical lemma. 
\begin{lemma}
\label{convex}
  Let $f$ be a convex function on $[0,1]$, say.
  Consider the problem
  \begin{equation}
    \label{eq:max-convex}
    f^* = \max \sum_{j = 1}^k f(x_j), \quad \text{subject to } x_j \ge 0 
\text{ and } \sum_{j = 1}^k x_j = 1.  
  \end{equation}
  Then the maximum value $f^*$ is obtained by allocating one $x_j$ to
  1 and all the others to 0, i.e. $f^* = (k-1) f(0) + f(1)$.
\end{lemma}
\begin{proof} For each $x_j$, $0 \le x_j \le 1$, the convexity of $f$
implies
\[
f(x_j) \le (1 - x_j) f(0) + x_j f(1).
\]
Summing this inequality over all indices gives 
\[
\sum_{j = 1}^k f(x_j) \le (k - 1) f(0) + f(1), 
\]
which is what we sought to establish. 
\end{proof}

\begin{corollary}
  Suppose that $f = \log F$ is a convex function on
  $[0,1]$, say, and consider
\begin{equation}
    \label{eq:max-convex2}
F^* = \max \prod_{j = 1}^k F(x_j), 
\quad \text{subject to } x_j \ge 0 \text{ and } \sum_{j = 1}^k x_j = 1.  
  \end{equation}
  Then the maximum value $F^*$ is obtained by allocating one $x_j$ to
  1 and all the others to 0, i.e. $F^* = (F(0))^{k-1} F(1)$.
\end{corollary}
\begin{proof}
Take the logarithm of $\prod_{j = 1}^k F(x_j)$ and apply 
Lemma \ref{convex}. 
\end{proof}

Note that both the lemma and the corollary hold for 'discrete'
functions; that is, suppose that $f(j)$ obeys
\begin{equation}
  \label{eq:discreteconvex}
f(j+1) - f(j) \ge f(j) - f(j-1), \quad j = 0, 1, 2, \ldots. 
\end{equation}
Then the maximum value of $\sum_{j = 1}^k f(n_j)$ where the $n_j$'s
are now integer values obeying $n_j \ge 0$ and $\sum_{j = 1}^k n_j =
n$ is of course achieved by taking all the $n_j$'s equal to zero but
one equal to $n$.

With these preliminaries in place, recall now the bound obtained in
Lemma \ref{fn-bound}, 
$$
F_n(\tau) \le G_{\tau/(1-\tau)}(n)
$$
where 
\begin{equation}
  \label{eq:G}
  G_u(n) = \begin{cases} u, & \text{$\log u\le 1 - n$},\\ 
 \exp((n-1) (\log(n-1) - \log \log (1/u) - 1)),    
& \text{$\log u >  1 - n$}. \end{cases} . 
\end{equation}
Note that we voluntarily exchanged the subscripts, namely, $\tau$ and
$n$ to reflect the idea that we shall view $G$ as a function of $n$
while $\tau$ will serve as a parameter. It is clear that $\log G$ is
convex and, therefore,
$$
G_u^* = \max \prod_{j = 1}^k G_u(n_j), 
\quad \text{subject to } n_j \ge 2 \text{ and } \sum_{j = 1}^k n_j = 2n
$$
obeys 
\begin{equation*}
   G_u^* = (G_u(2))^{k-1} G_u(2n-2k+2).
\end{equation*}

Set $G = G_{\tau/(1-\tau)}$ for short. Then for any equivalence class
such that $|A /\sim| = k$, the above argument yields
$$
\prod_{A' \in A/\sim} F_{|A'|}(\tau) \le \prod_{A' \in A/\sim}
G(|A'|) \le [G(2)]^{k-1} \, G(2n-2k+2), 
$$
which, on the one hand, gives 
\begin{equation*}
 \E\tr(H^{2n}_0)  \le \sum_{k = 1}^{n} N^{k} \, |T|^{2n - k + 1} \, P(2n,k) \, 
  [G(2)]^{k-1} \, G(2n-2k+2). 
\end{equation*}
On the other hand, $P(2n,k) \le \phi^{2n}\, (2n-1) \ldots
(2n-2k+1)$ (see Lemma \ref{Pnk}) and, therefore,
\begin{equation}
  \label{eq:expected5}
   \E\tr(H^{2n}_0)  \le |T| \, \lambda^{2n} \, \sum_{k = 1}^n f(k), 
\end{equation}
where
\begin{equation}
  \label{eq:fk}
  f(k) := N^k \, |T|^{2n-k} \,[(2n-1) \ldots (2n-2k+1)] 
\, [G(2)]^{k-1} \, G(2n-2k+2)  
\end{equation}

We prove that the summand $f$ is in some sense convex.
\begin{lemma}
\label{fdiscreteconvex}
   For each $k \le n-1$, $f$ obeys 
$$
f(k+1) - f(k) \ge f(k) - f(k-1). 
$$
\end{lemma}
As a consequence of this lemma, the maximum of $f(k)$, $1\le k \le n$
is of course attained at either the left-end point ($k = 1$) or
the right-end point ($k = n$); in short, 
\[
f(k) \le \max(f(1), f(n)), \quad \forall 1 \le k \le n. 
\]

\begin{proof}
We need to establish that for each $1 \le k \le
n-1$, 
$$
\frac{f(k+1)}{f(k)} + \frac{f(k-1)}{f(k)} \ge 2. 
$$
Observe that 
$$
\frac{f(k+1)}{f(k)} = \alpha \rho_{k+1}, \quad \frac{f(k)}{f(k-1)}
= \alpha^{-1} \rho_{k-1},
$$
with $\alpha = N\, G(2)/|T|$ and 
$$
\rho_{k+1} = (2n-2k-1) \cdot \frac{G(2n-2k)}{G(2n-2k+2)}, \quad
\rho_{k-1} = \frac{1}{(2n-2k+1)} \cdot \frac{G(2n-2k+4)}{G(2n-2k+2)}.
$$
Clearly
$$
\alpha \rho_{k+1} + \alpha^{-1} \rho_{k-1} \ge 2\, \sqrt{
  \rho_{k+1} \, \rho_{k-1}}
$$
and, therefore, it is sufficient to establish that $\rho_{k+1} \,
\rho_{k-1} \ge 1$. Put $m = n - k$, then
\begin{eqnarray*}
\rho_{k+1} \, \rho_{k-1} & = & \frac{m-1}{m+1} \cdot \frac{G(m)
\, G(m+4)}{[G(m+2)]^2} \\
& = & \frac{m-1}{m+1} \cdot 
\frac{(m-1)^{m-1} \,
  (m+3)^{m+3}}{(m+1)^{2(m+1)}}. 
\end{eqnarray*}
It is now a simple exercise to check that for each $m \ge 2$, the
logarithm of the right-hand is nonnegative, i.e.
$$
m \log (m-1) + (m+3) \log(m+3) - (2m + 3) \log(m+1) \ge 0.  
$$
We omit the proof of this fact. 
\end{proof}

\subsection{Proof of Theorem \ref{moment-bound}}

The previous section established 
$$
\E\tr(H^{2n}_0) \le |T| \, \phi^{2n} \cdot n \cdot
\max(f(1),f(n))
$$
where letting $c_\tau := e \log((1-\tau)/\tau)$
$$
f(1) = (2n-1) \, G(2n) \, N \, |T|^{2n-1} = (2n-1)^{2n} \,
c_\tau^{-(2n-1)}\, N \, |T|^{2n-1}.
$$
and
$$
f(n) = [(2n-1) \times (2n-3) \ldots \times 1] \cdot
\left(\frac{\tau}{1-\tau}\right)^n \, N^n \, |T|^n.
$$
This is exactly the content of Theorem \ref{moment-bound}. 


\section{Numerical Experiments}
\label{sec:numexp}

In this section, we present numerical experiments in order to derive
empirical bounds on $|\supf|$ relative to $|\Omega|$ for a signal $f$
supported on $\supf$ to be the unique minimizer of $(P_1)$.  The
results can be viewed as a set of practical guidelines for situations
where one can expect perfect recovery from partial Fourier information
using convex optimization.

\begin{figure}
\centerline{
\begin{tabular}{cc}
\raisebox{1in}{\rotatebox{90}{$|\Omega|$}}
\includegraphics[width=3in]{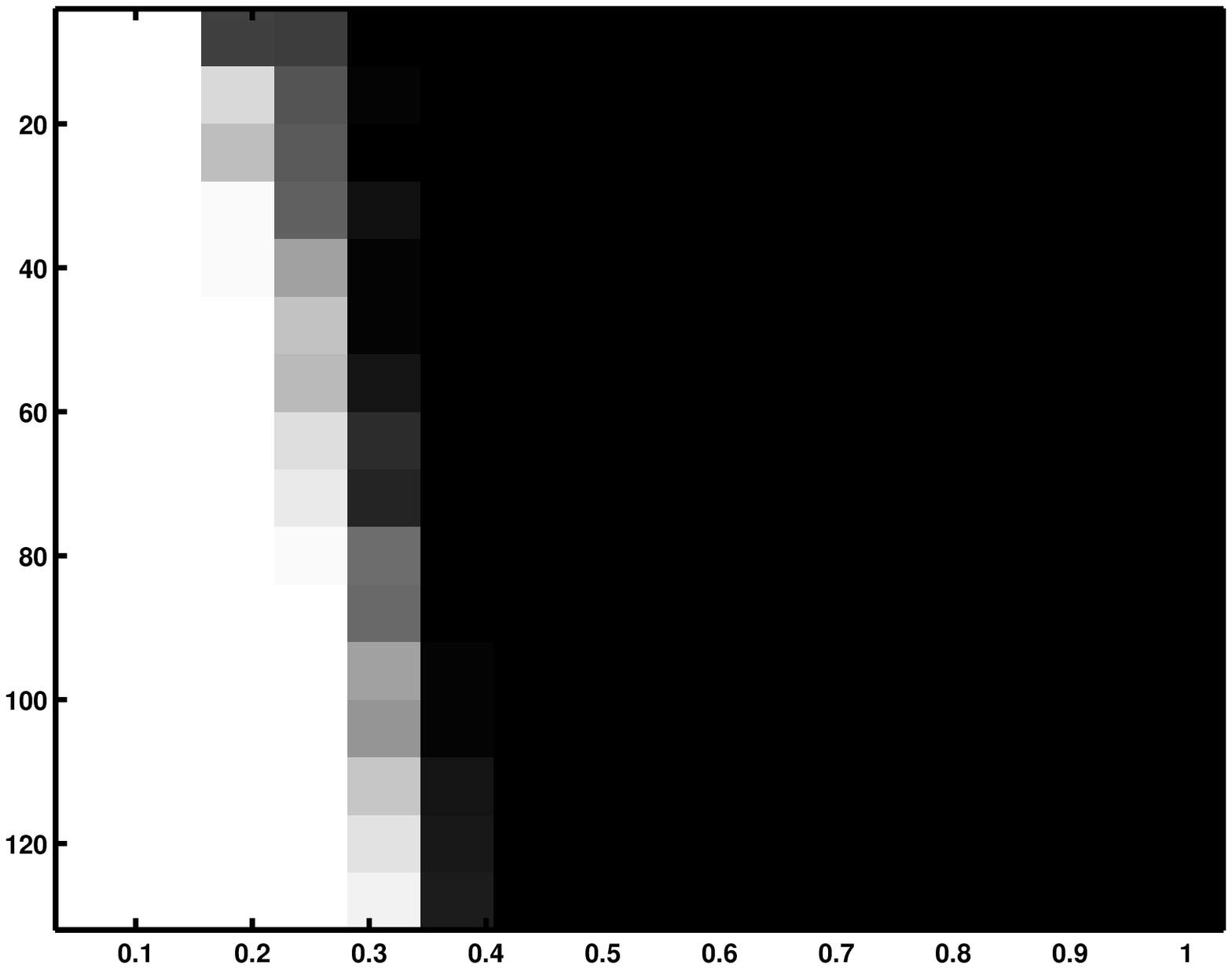} &
\raisebox{.8in}{\rotatebox{90}{\sf recovery rate}}
\includegraphics[width=3in]{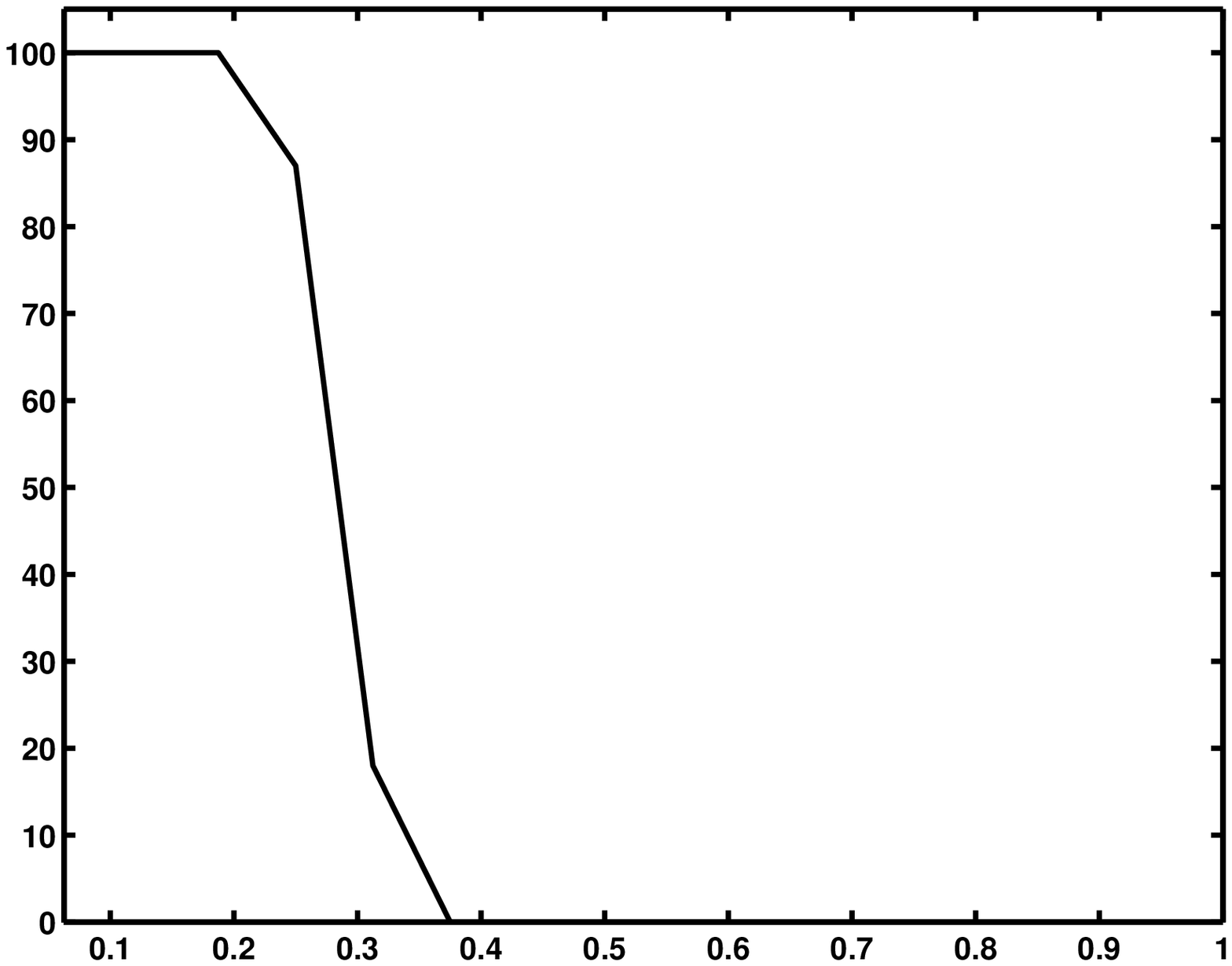} \\
$|\supf|/|\Omega|$ & $|\supf|/|\Omega|$ \\[3mm]
(a) & (b)
\end{tabular}
}
\caption{Recovery experiment for $N=512$.  (a) The image intensity
  represents the percentage of the time solving $(P_1)$ recovered the
  signal $f$ exactly as a function of $|\Omega|$ (vertical axis) and
  $|\supf|/|\Omega|$ (horizontal axis); in white regions, the signal
  is recovered approximately $100\%$ of the time, in black regions,
  the signal is never recovered.  For each $|\supf|,|\Omega|$ pair,
  $100$ experiments were run.
(b) Cross-section of the image in (a) at $|\Omega| = 64$.  We can see
that we have perfect recovery with very high probability for
$|\supf|\leq 16$.  }
\label{fig:recover512}
\end{figure}

Our experiments are of the following form:
\begin{enumerate}
\item Choose constants $N$ (the length of the signal), $N_{t}$
  (the number of spikes in the signal), and $N_{\omega}$ (the number
  of observed frequencies).
\item Randomly generate the subdomain $\supf$ by sampling
  $\{0,\ldots,N-1\}$ $N_{t}$ times without replacement (we have
  $|\supf| = N_{t}$).
\item Randomly generate $f$ by setting $f(t)=0, t\in\supf^c$ and
  drawing both the real and imaginary parts of $f(t), t\in\supf$ from
  independent Gaussian distributions with mean zero and variance
  one\footnote{The results here, as in the rest of the paper, seem to
    rely only on the sets $\supf$ and $\Omega$.  The actual values
    that $f$ takes on $\supf$ can be arbitrary; choosing them to be
    random emphasizes this.  Figures~\ref{fig:recover512} remain the
    same if we take $f(t)=1,t\in\supf$, say.}.
\item Randomly generate the subdomain $\Omega$ of observed frequencies
  by again sampling $\{0,\ldots,N-1\}$ $N_{\omega}$ times without
  replacement ($|\Omega| = N_{\omega}$).
\item Solve $(P_1)$, and compare the solution to $f$.
\end{enumerate}

The $\ell_1$-norm is not strictly convex, so solving $(P_1)$ using a
Newton-type method that relies on local quadratic approximations of
$\|\cdot\|_{\ell_1}$ is problematic.  Instead, we use a very simple
gradient descent with projection algorithm.  The number of iterations
needed for convergence is high (on the order of $10^5$), but since we
can rapidly project onto the constraint set (using two fast Fourier
transforms), each iteration takes a short amount of time.  As an
indication, the algorithm typically converges in less than $10$
seconds on a standard desktop computer for signals of length $N=1024$.

Figure~\ref{fig:recover512} illustrates the recovery rate for varying
values of $|\supf|$ and $|\Omega|$ for $N=512$.  From the plot, we
can see that for $|\Omega|\geq 32$, if $|\supf| \leq |\Omega|/5$, we
recover $f$ perfectly about $80\%$ of the time.  For $|\supf| \leq
|\Omega|/8$, the recovery rate is practically $100\%$. We remark that
these numerical results are consistent with earlier findings
\cite{Loh}.

One source of slack in the theoretical analysis is the way in which we
choose the polynomial $P(t)$ (as in \eqref{P-def}).
Theorem~\ref{duality} states that $f$ is a minimizer of $(P_1)$ if and
only if there exists {\em any} trigonometric polynomial that has
$P(t)=\sgn(f)(t), t\in\supf$ and $|P(t)|<1, t\in\supf^c$.  In
\eqref{P-def} we choose $P(t)$ that minimizes the $\ell_2$ norm on
$\supf^c$ under the linear constraints $P(t)=\sgn(f)(t),t\in\supf$.
However, the condition $|P(t)| < 1$ suggests that a minimal
$\ell_\infty$ choice would be more appropriate (but is seemingly
intractable analytically).

\begin{figure}
\centerline{
\begin{tabular}{cc}
\raisebox{1in}{\rotatebox{90}{$|\Omega|$}}
\includegraphics[width=3in]{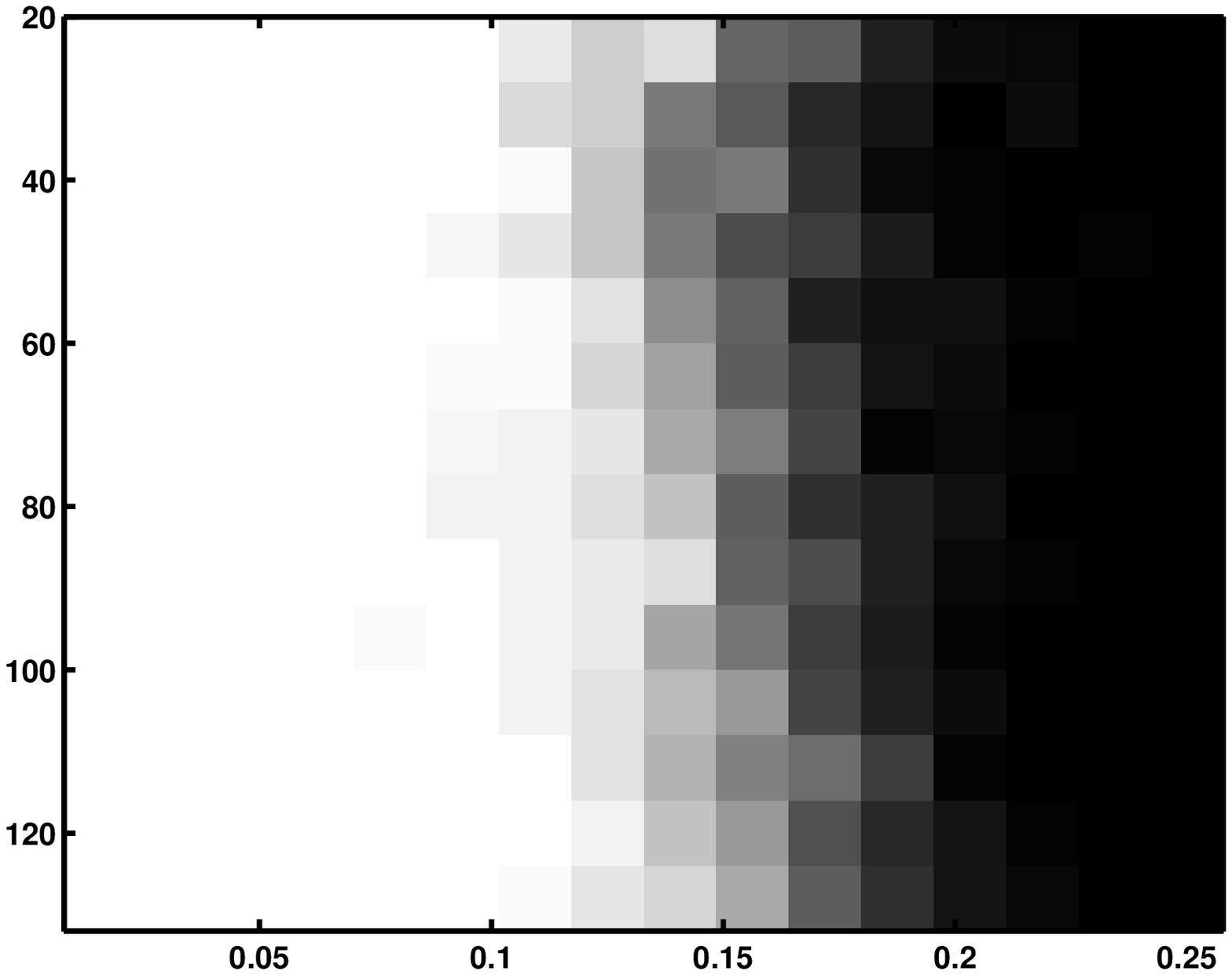} &
\raisebox{.8in}{\rotatebox{90}{\sf \% suff. true}}
\includegraphics[width=3in]{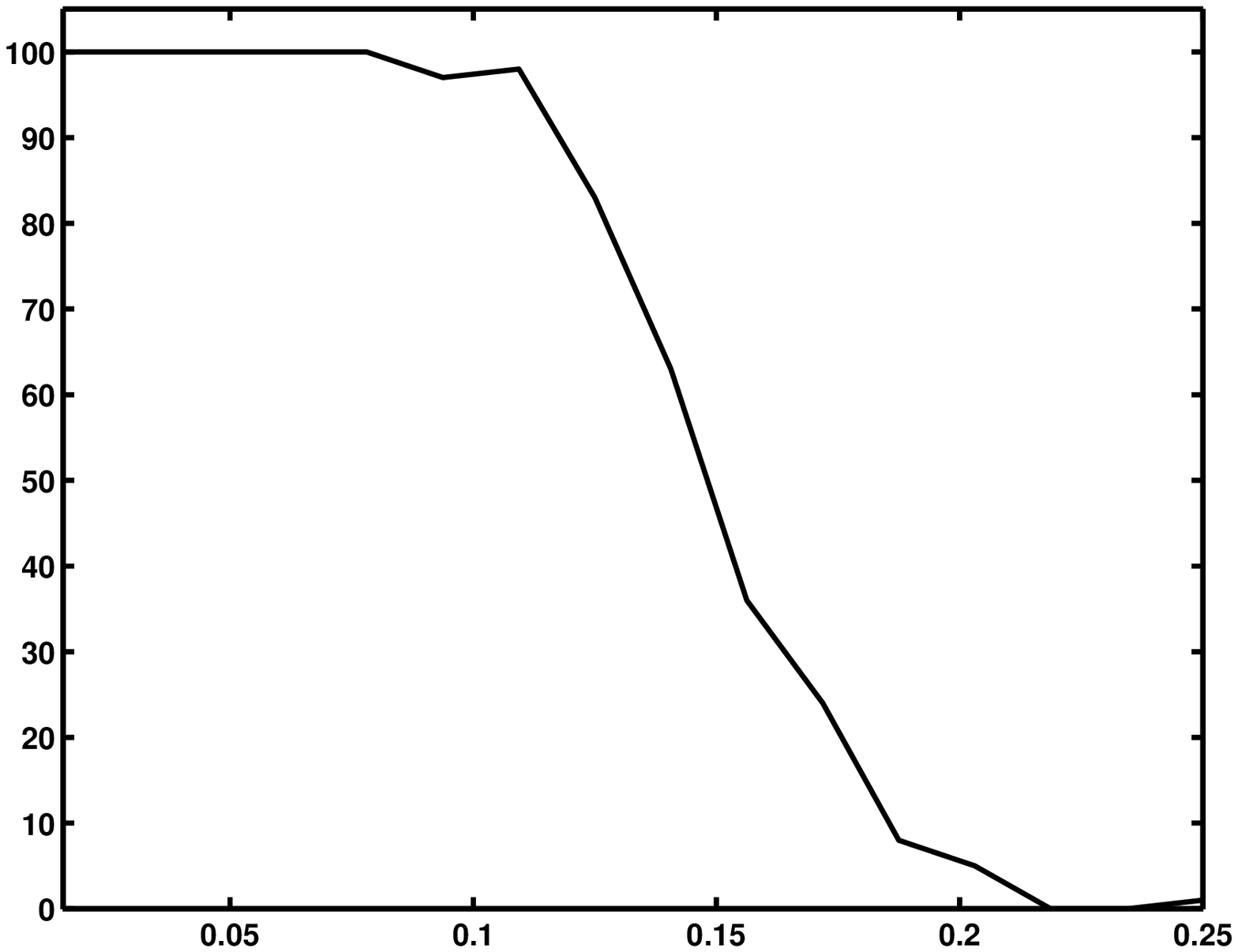} \\
$|\supf|/|\Omega|$ & $|\supf|/|\Omega|$ \\[3mm]
(a) & (b)
\end{tabular}
}
\caption{Sufficient condition test for $N=512$. (a) The image
  intensity represents the percentage of the time $P(t)$ chosen as in
  \eqref{P-def} meets the condition $|P(t)|<1, t\in\supf^c$. (b) A
  cross-section of the image in (a) at $|\Omega| = 64$.  Note that the
  axes are scaled differently than in Figure~\ref{fig:recover512}.  }
\label{fig:polysuff}
\end{figure}

Figure~\ref{fig:polysuff} illustrates how often the sufficient
condition of $P(t)$ chosen as \eqref{P-def} meets the constraint
$|P(t)|<1, t\in\supf^c$ for the same values of $\tau$ and $|\supf|$.
The empirical bound on $\supf$ is stronger by about a factor of two;
for $|\supf| \leq |\Omega|/10$, the success rate is very close to
$100\%$.

\begin{figure}
\centerline{
\begin{tabular}{ccc}
\includegraphics[width=2in]{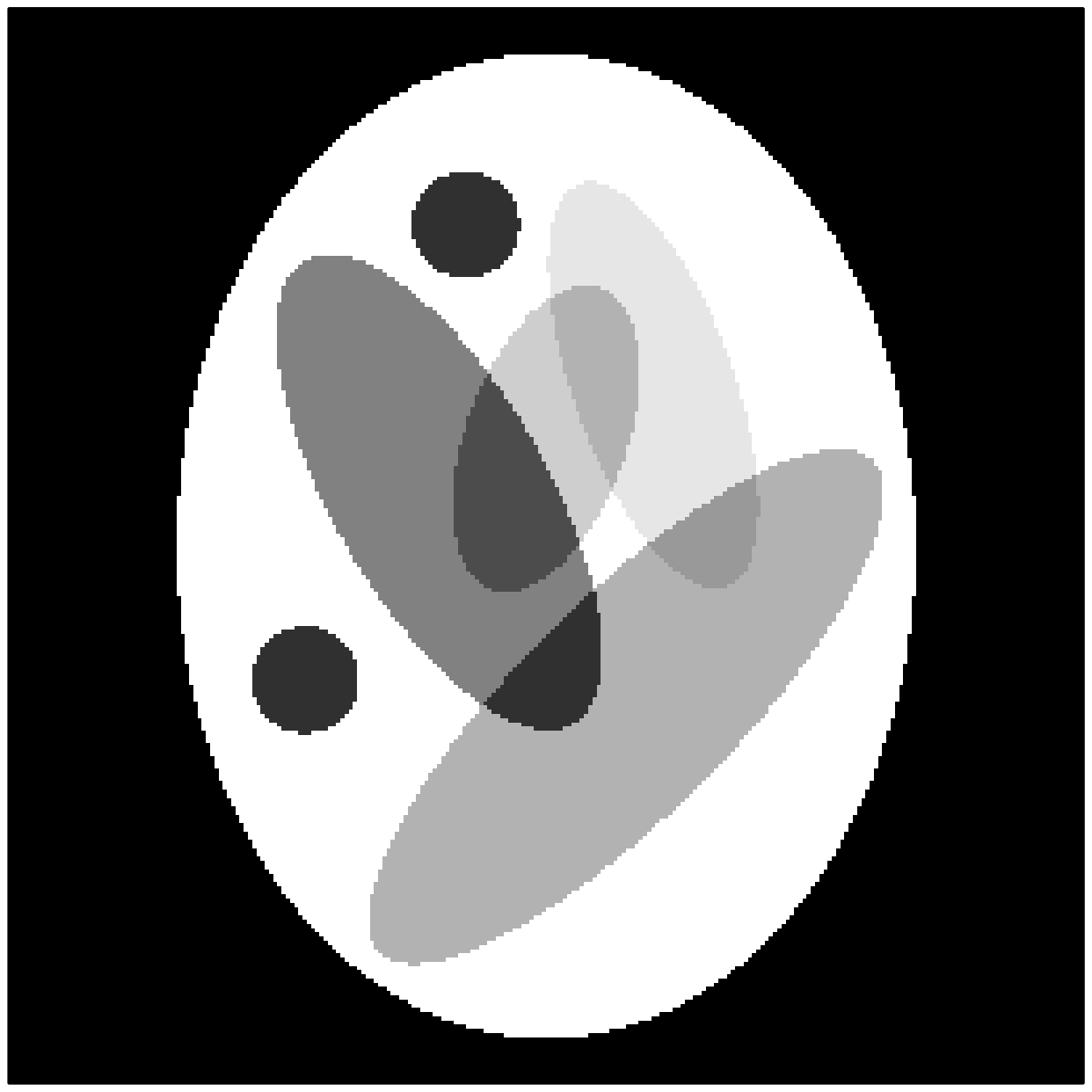} &
\includegraphics[width=2in]{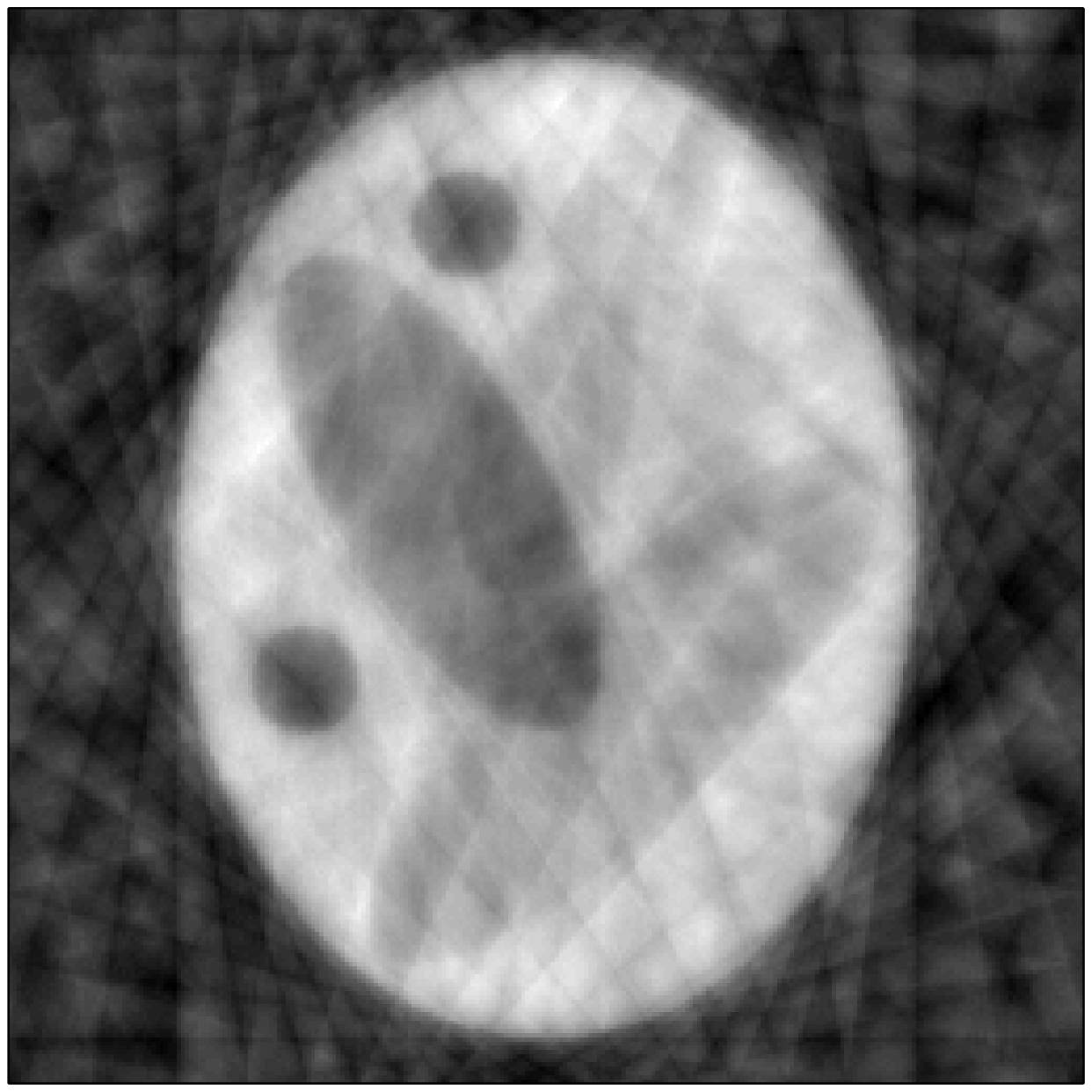} &
\includegraphics[width=2in]{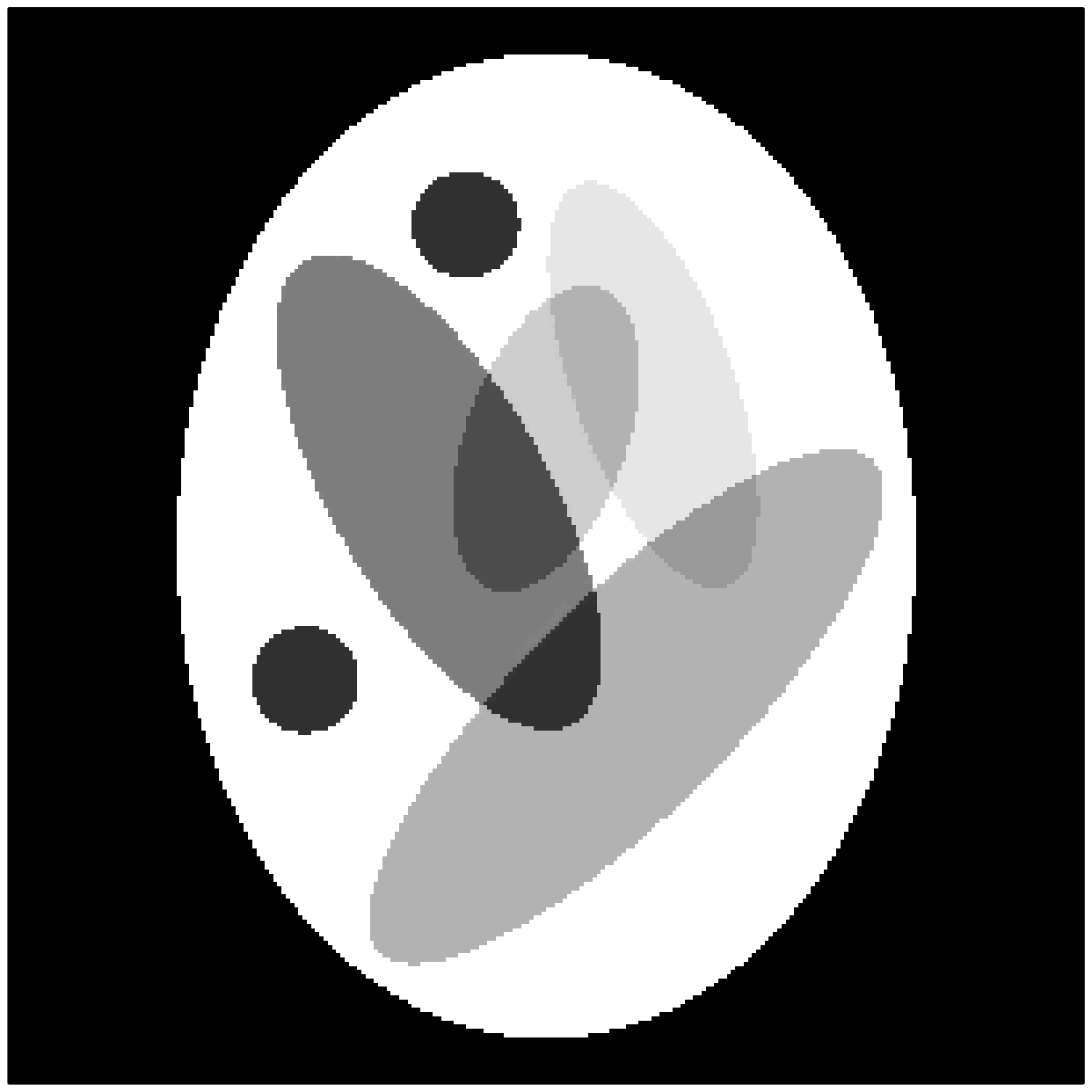} \\
(a) & (b) & (c) \\[3mm]
\includegraphics[width=2in]{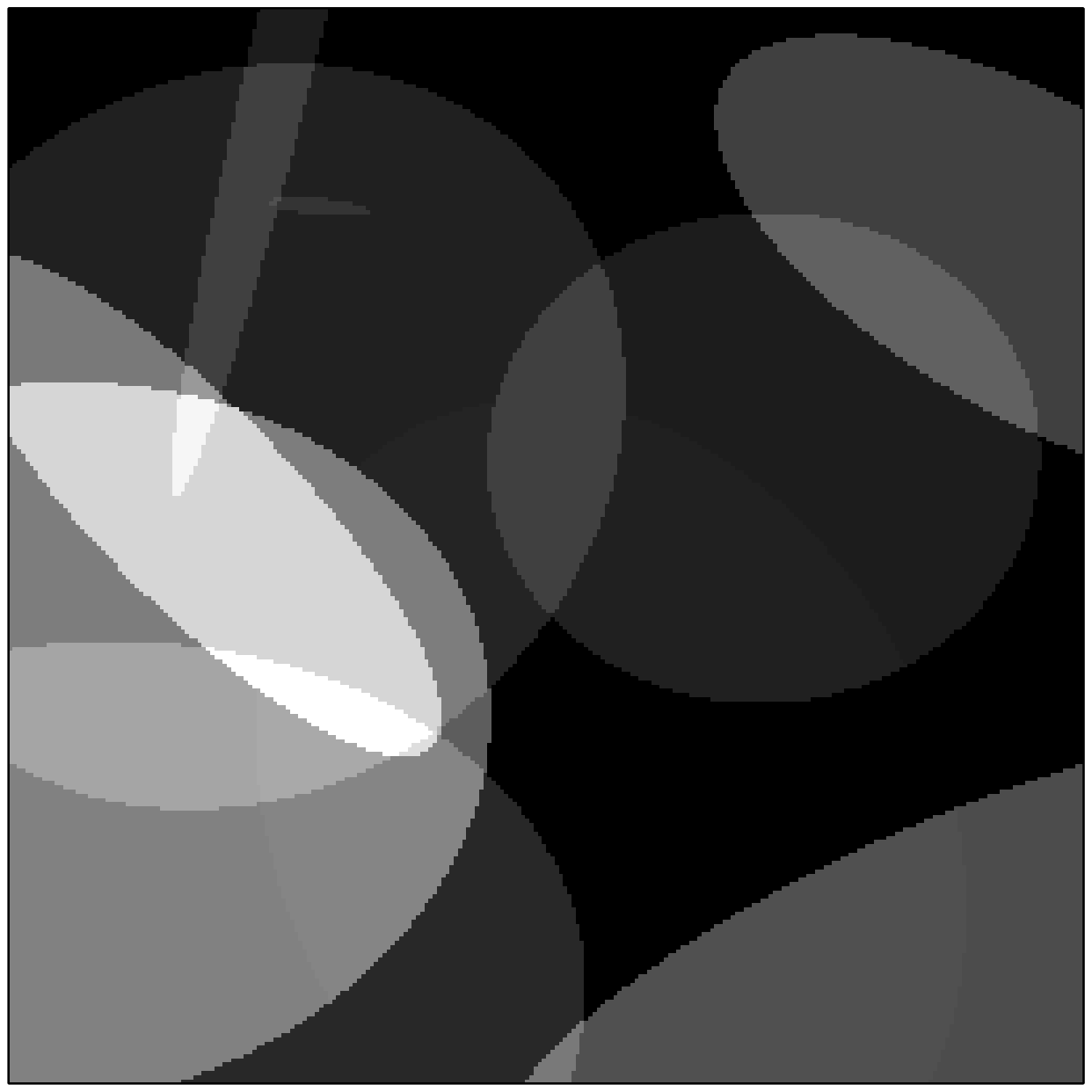} &
\includegraphics[width=2in]{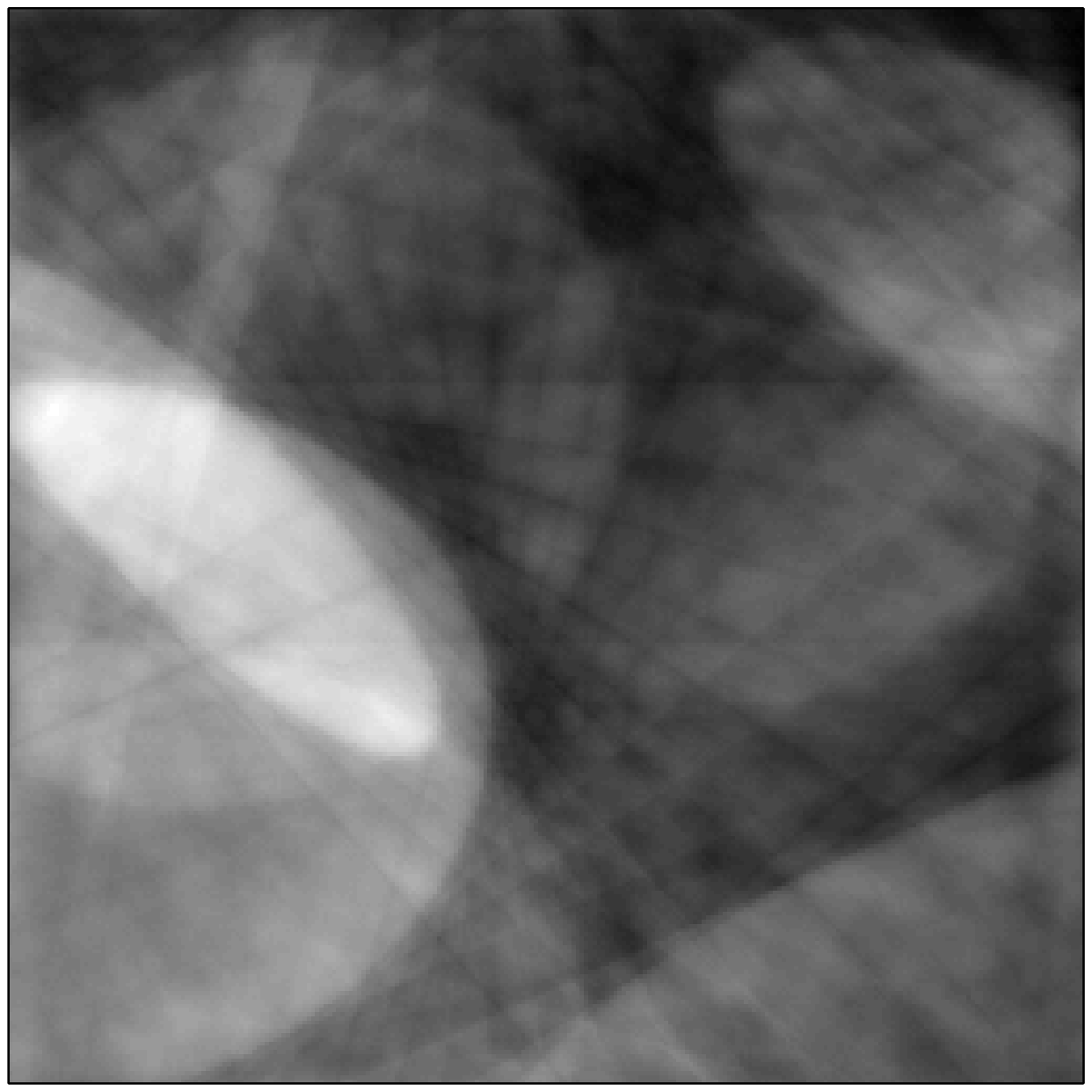} &
\includegraphics[width=2in]{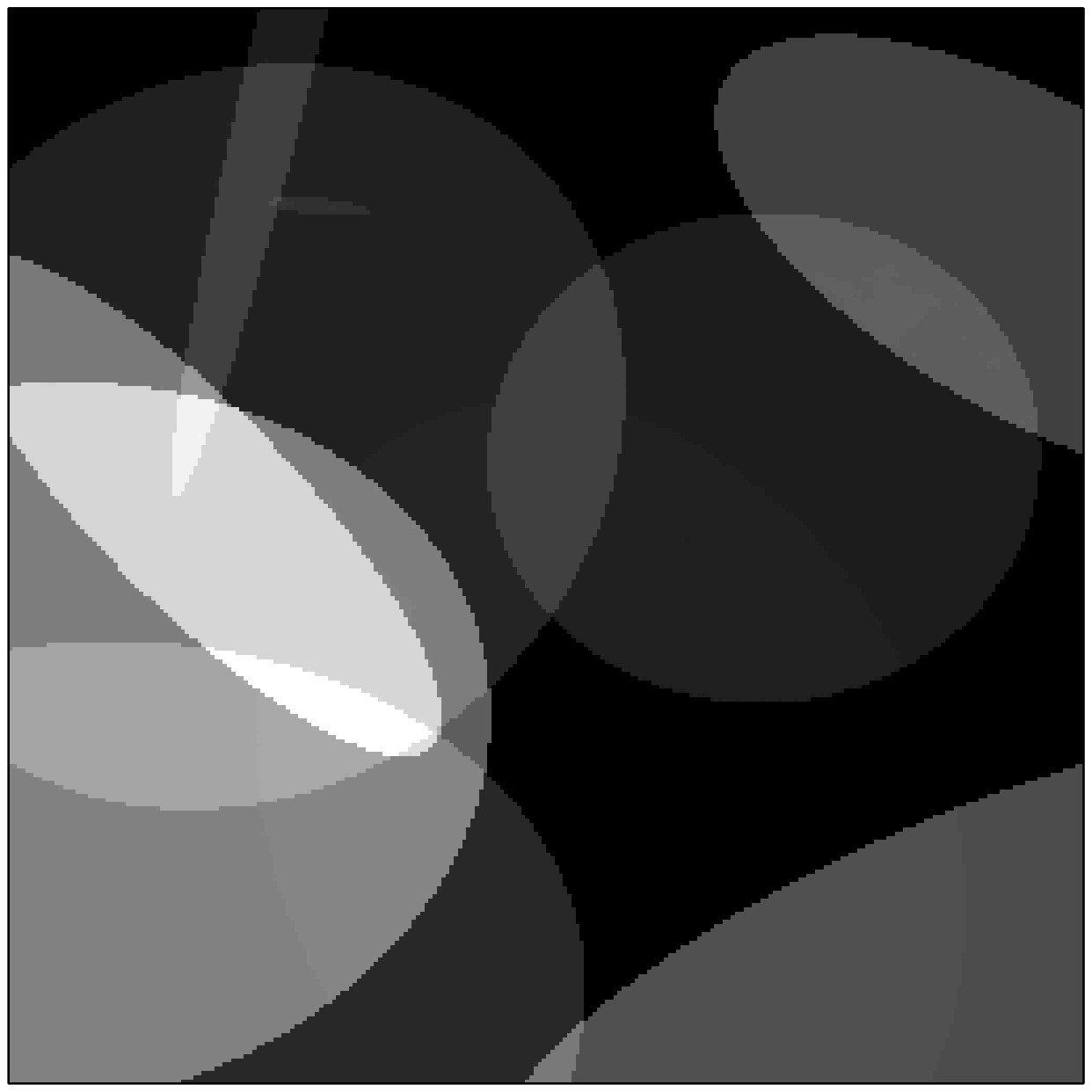} \\
(d) & (e) & (f) 
\end{tabular}
}
\caption{Two more phantom examples for the recovery problem discussed in Section~\ref{sec:puzphantom}.  On the left is the original phantom ((d) was created by drawing ten ellipses at random), in the center is the minimum energy reconstruction, and on the right is the minimum total-variation reconstruction.  The minimum total-variation reconstructions are exact.}
\label{fig:morephantoms}
\end{figure}

As a final example of the effectiveness of this recovery framework, we
show two more results of the type presented in
Section~\ref{sec:puzphantom}; piecewise constant phantoms
reconstructed from Fourier samples on a star.  The phantoms, along
with the minimum energy and minimum total-variation reconstructions
(which are exact), are shown in Figure~\ref{fig:morephantoms}.  Note
that the total-variation reconstruction is able to recover very subtle
image features; for example, both the short and skinny ellipse in the
upper right hand corner of Figure~\ref{fig:morephantoms}(d) and the
very faint ellipse in the bottom center are preserved. (We invite the
reader to check \cite{BreslerDelaney} for related types of
experiments.)

\section{Discussion}
\label{sec:discussion}

We would like to close this paper by offering a few comments about the
results obtained in this paper and by discussing the possibility of
generalizations and extensions. 

\subsection{Stability}

In the introduction section, we argued that even if one knew the
support $T$ of $f$, the reconstruction might be unstable. Indeed with
knowledge of $T$, a reasonable strategy might be to recover $f$ by the
method of least-squares, namely,
\[
f = ({\cal F}_{T \goto \Omega}^* {\cal F}_{T \goto \Omega})^{-1} \, 
{\cal F}_{T \goto \Omega}^* \, \hat f|_\Omega. 
\]
In practice, the matrix inversion might be problematic. Now observe that
with the notations of this paper
\[
{\cal F}_{T \goto \Omega}^* {\cal F}_{T \goto \Omega} \propto
I_T - \frac{1}{|\Omega|} H_0. 
\]
Hence, for stability we would need $\frac{1}{|\Omega|} H_0 \le 1 -
\delta$ for some $\delta > 0$. This is of course exactly the problem
we studied, compare Theorem \ref{super-useful}. In fact, selecting
$\alpha_M$ as suggested in the proof of our main theorem (see section
\ref{sec:proof}) gives $\frac{1}{|\Omega|} H_0 \le .42$ with
probability at least $1 - O(N^{-M})$. This shows that selecting $|T|$
as to obey \eqref{supp}, $|T| \approx |\Omega|/\log N$ actually
provides stability.

\subsection{Robustness}

An important question concerns the robustness of the reconstruction
procedure vis a vis measurement errors.  For example, we might want to
consider the model problem which says that instead of observing the
Fourier coefficients of $f$, one is given those of $f + h$ where $h$
is some small perturbation. Then one might still want to reconstruct
$f$ via
\[
f^\sharp = \text{argmin } \|g\|_{\ell_1}, \qquad \hat g(\omega) = \hat
f(\omega) + \hat h(\omega), \quad \forall \omega \in \Omega.
\]
In this setup, of course, one cannot expect exact recovery. Instead,
one would like to know whether or not our reconstruction strategy is
well-behaved or more precisely, how far is the minimizer $f^\sharp$
from the true object $f$. In short, what is the typical size of the
error?  Our preliminary calculations suggest that the reconstruction
is robust in the sense that the error $ \|f - f^\sharp\|_1$ is small
for small perturbations $h$ obeying $\|h\|_1 \le \delta$, say. We hope
to be able to report on these early findings in a follow-up paper.

\subsection{Extensions}

Finally, work in progress shows that similar exact reconstruction
phenomena hold for other synthesis/measurement pairs.  Suppose one is
given a pair of of bases $({\cal B}_1, {\cal B}_2)$ and randomly
selected coefficients of an object $f$ in one basis, say ${\cal B}_2$.
(From this broader viewpoint, the special cases discussed in this
paper assume that ${\cal B}_1$ is the canonical basis of $\RR^N$ or 
$\RR^N \times \RR^N$ (spikes in 1D, 2D), or is the basis of Heavysides as
in the Total-variation reconstructions, and ${\cal B}_2$ is the
standard 1D, 2D Fourier basis.) Then, it seems that $f$ can be
recovered exactly provided that it may be synthesized as a sparse
superposition of elements in ${\cal B}_1$.  The relationship between
the number of nonzero terms in ${\cal B}_1$ and the number of observed
coefficients depends upon the {\em incoherence} between the two bases
\cite{DonohoHuo}.  The more incoherent, the fewer coefficients
needed. Again, we hope to report on such extensions in a separate
publication.

\section{Appendix}

\subsection{Proof of Theorem \ref{super-useful}}

\begin{figure}
  \centering
\begin{tabular}{cc}
\includegraphics[width=2.5in]{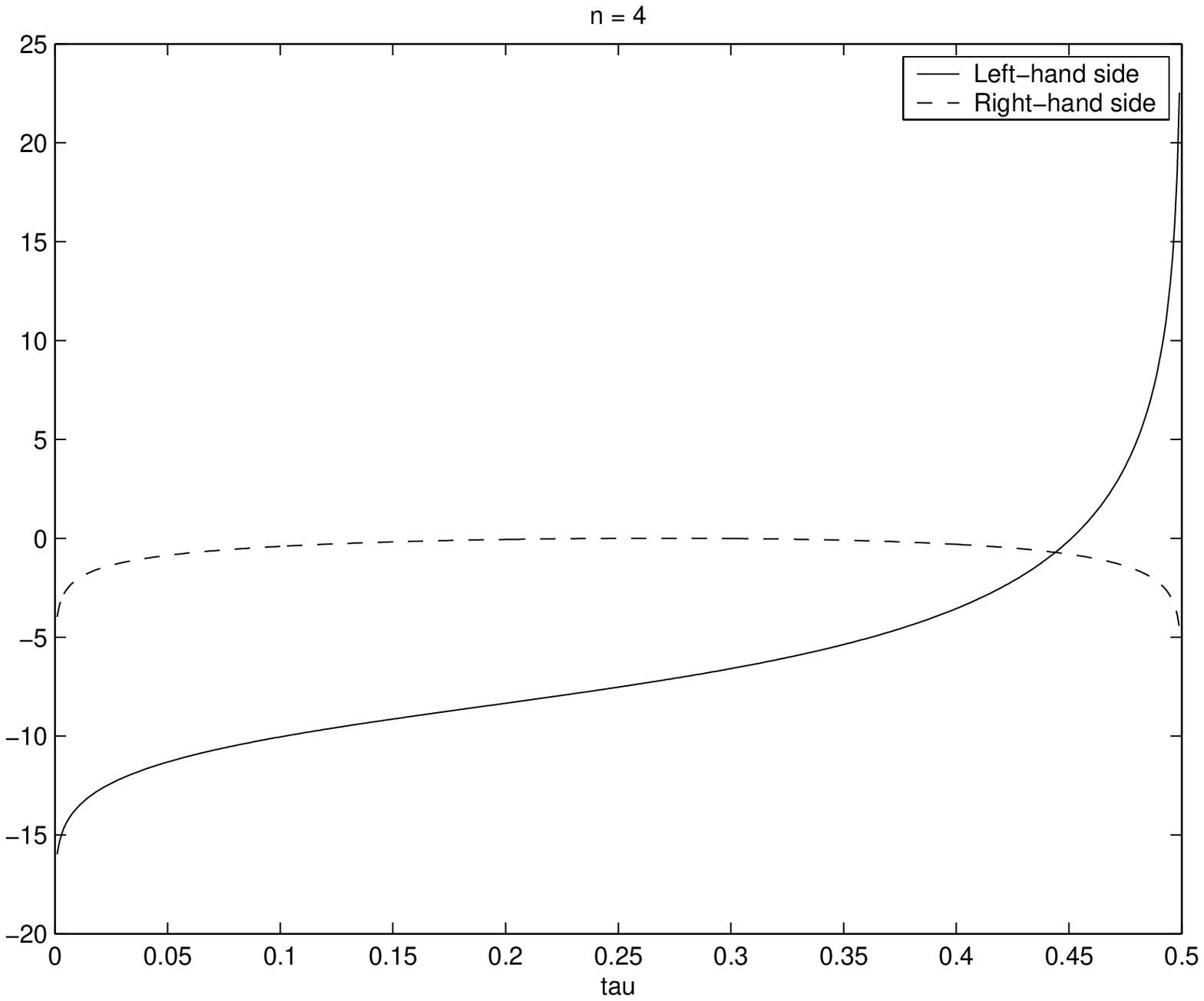} &
\includegraphics[width=2.5in]{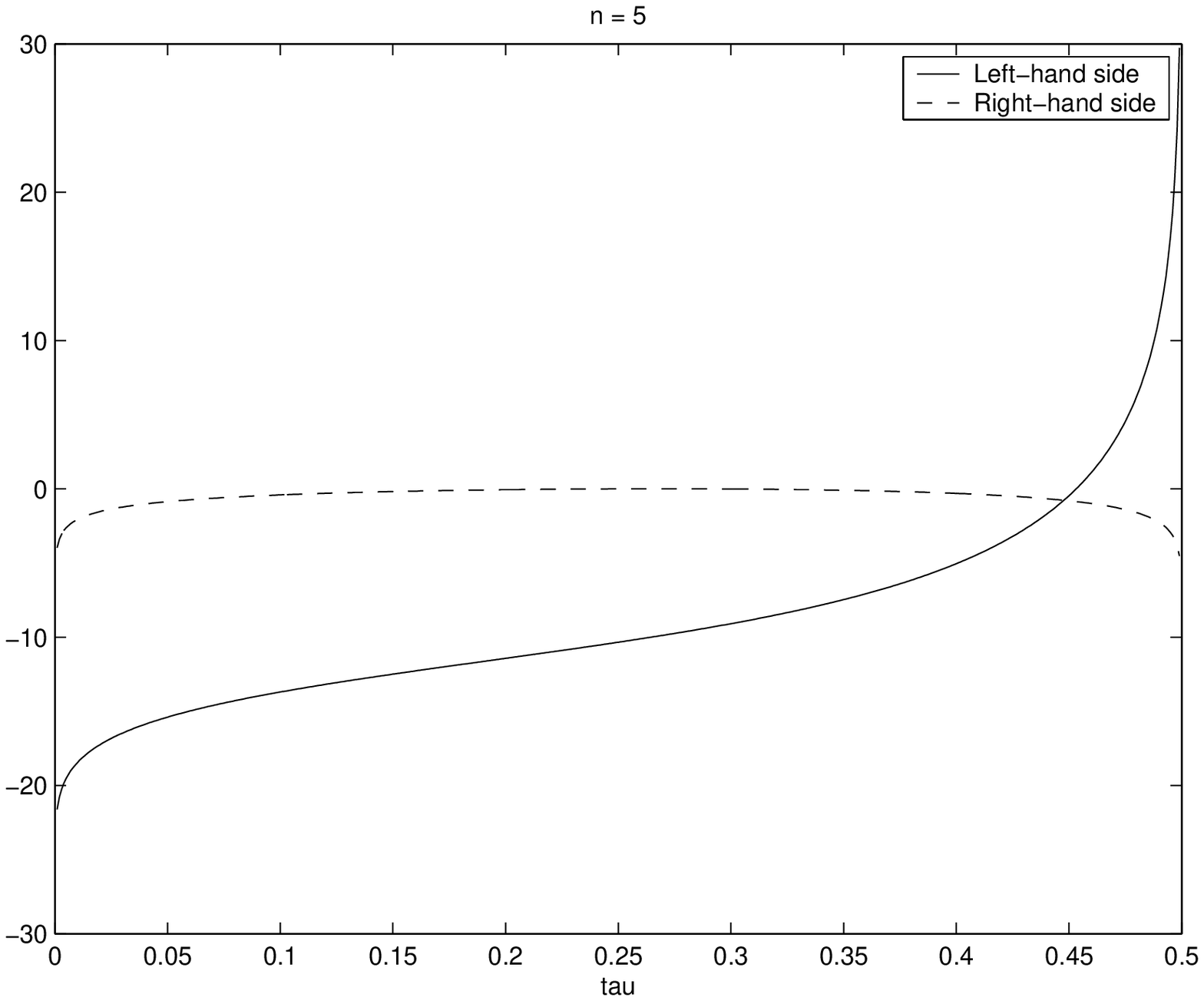} \\
$n = 4$ & $n = 5$\\
\end{tabular}
\caption{Behavior of the left and right-hand side of 
\eqref{eq:log-ineq} for two values of $n$}
\label{fig:plot} 
\end{figure}

We need to prove that for $\tau \le .44$ and $n \ge 4$,
$$
(2n-1)^{2n} \, (c_\tau)^{-(2n-1)} \, N \, |T|^{2n-1} \le n^{n} \,
2^{n+1} \, e^{-n} \, \left(\frac{\tau}{1-\tau}\right)^n \, N^n \,
|T|^n;
$$
Now $(2n-1)^{2n} = (2n)^{2n} e^{-1} \epsilon_n$ where $\epsilon_n
\le e^{1/2n}$, say. We may then rewrite the previous inequality as 
$$
\epsilon_n \, (2e)^{n-1} \, n^n (c_\tau)^{-2(n-1)} \, 
(1-\tau)^{n-1} |T|^{n-1} \le |\Omega|^{n-1} \, s_\tau
$$
where $s_\tau = \frac{\tau}{1-\tau} c_\tau$. Because $|T| \le
\frac{\alpha^2}{\gamma^2} \frac{|\tau N|}{n}$, $0 < \alpha < 1$, it is
sufficient to check that 
$$
\epsilon_n \, n \, r_\tau^{n-1} \le s_\tau, \quad 
r_\tau = \frac{2 \alpha^2 e (1 - \tau)}{\gamma^2 \, c_\tau^2}.
$$
Note that plugging the value of $\gamma$ gives $r_\tau =
(1-\tau)^3/(e \, \alpha ^ 2 \, \phi^2
\left[\log(\frac{1-\tau}{\tau})\right]^2)$ (recall $\phi =
(1+\sqrt{5})/2$). In other words, we want
\begin{equation}
  \label{eq:log-ineq}
  (n - 1) \log r_\tau + \log n + \frac{1}{2n} \le \log s_\tau. 
\end{equation}
Figure~\ref{fig:plot} illustrates the behavior of both the
left-hand side and the right-hand side with $\alpha = 1$. 
Simple numerical calculations show that with $\alpha = 1$,
\eqref{eq:log-ineq} holds for $\tau \le .44$ and $n \ge 4$, as
claimed.

\subsection{Proof of Lemma \ref{second-moment}}

Set $e^{i\phi} = \sgn(f)$ for and fix $K$. Using \eqref{T-def}, we
have
$$
[(H \iota^*)^{n+1} e^{i\phi}](t_0) = \sum_{t_1, \ldots, t_{n+1} \in
  T: \, t_j \neq t_{j+1} \text{ for } j = 0, \ldots, n} \,\,\,
\sum_{\omega_0, \ldots, \omega_{n} \in \Omega} \,\,\, e^{i
  \sum_{j=0}^{n} \omega_j (t_j - t_{j+1})} \,\, e^{i\phi(t_{n+1})},
$$
and, for example, 
\begin{multline*}
  |[(H \iota^*)^{n+1} e^{i\phi}](t_0)|^2 = \sum_{t_1, \ldots, t_{n+1} \in T: \,
    t_j \neq t_{j+1} \text{ for } j = 0, \ldots, n \atop t'_1, \ldots,
    t'_{2n} \in T: \, t'_j \neq t'_{j+1} \text{ for } j = 0, \ldots,
    n} \,\, e^{i\phi(t_{n+1})} \,
  e^{-i\phi(t'_{n+1})}\\
  \sum_{\omega_0, \ldots, \omega_{n} \in \Omega \atop \omega'_0,
    \ldots, \omega'_{n} \in \Omega} \,\, e^{i \sum_{j=0}^{n}
    \omega_j(t_j - t_{j+1})} \, e^{-i \sum_{j=0}^{n} \omega'_j (t'_j -
    t'_{j+1})}.
\end{multline*}
One can calculate the $2K$th moment in a similar fashion. Put $m: = K(n+1)$ 
and 
$$
\bomega := (\omega^{(k)}_j)_{k,j}, 
\quad \bt = (t^{(k)}_j)_{k,j} \in T^{2K(n+1)}, 
\quad 1 \leq j \leq n+1 \text{ and } 1 \leq k \leq 2K 
$$
With these notations, we have 
$$
  |[(H \iota^*)^{n+1} g](t_0)|^{2K} = 
\sum_{\bt \in T^{2m}: t^{(k)}_j \neq t^{(k)}_{j+1}} 
  \sum_{\bomega \in \Omega^{2m}} \,\, e^{i \sum_{k = 1}^{2K} (-1)^{k}
    \phi(t^{(k)}_{n+1})} \,\, e^{i \sum_{k = 1}^{2K} \sum_{j=0}^{n}
    (-1)^{k} \omega^{(k)}_j(t^{(k)}_j - t^{(k)}_{j+1})}, 
  $$
  where we adopted the convention that $x^{(k)}_0 = x_0$ for all $1
  \leq k \leq 2K$ and where it is understood that the condition
  $t^{(k)}_j \neq t^{(k)}_{j+1}$ is valid for $0 \leq j \leq n$. 
  
  Now the calculation of the expectation goes exactly as in section
  \ref{sec:moments}. Indeed, we define an equivalence relation
  $\sim_{\bomega}$ on the finite set $A := \{ 0, \ldots, n \}
  \times \{1, \ldots, 2K\}$ by setting $(j,k) \sim (j',k')$ if
  $\omega_j^{(k)} = \omega_{j'}^{(k')}$ and observe as before that
  $$
  \E\left[\prod_{j,k} I_{\omega_j^{(k)}}\right] = \tau^{|A/\sim|};$$
  that is, $\tau$ raised at the power that equals the number of
  distinct $\omega$'s and, therefore, we can write the expected value
  $m(n;K)$ as
\begin{multline*}
  m(n;K) = \sum_{\bt \in T^{2m}: t^{(k)}_j \neq t^{(k)}_{j+1}} \,\,
\,\, e^{i \sum_{k = 1}^{2K}
    (-1)^{k} \phi(t^{(k)}_{n+1})} \,  \sum_{\sim \in {\cal P}(A)} 
 \,\,  \tau^{|A/\sim|}\\
  \sum_{\bomega \in \Omega^(\sim)} \,\,\,  e^{i \sum_{k = 1}^{2K}
    \sum_{j=0}^{n} (-1)^{k} \omega^{(k)}_j(t^{(k)}_j - t^{(k)}_{j+1})}
\end{multline*}
As before, we follow Lemma \ref{expected1} and rearrange this as
\begin{multline*}
  m(n;K) =  \sum_{\sim \in {\cal P}(A)} \,\, 
\sum_{\bt \in T^{2m}: t^{(k)}_j \neq t^{(k)}_{j+1}} 
\,\, e^{i \sum_{k = 1}^{2K}
    (-1)^{k} \phi(t^{(k)}_{n+1})} \,
\prod_{A' \in A/\sim}  F_{|A'|}(\tau) \\
\sum_{\bomega \in \Omega^(\sim)} \,\, e^{i \sum_{k = 1}^{2K} \sum_{j=0}^{n}
    (-1)^{k} \omega^{(k)}_j(t^{(k)}_j - t^{(k)}_{j+1})}
\end{multline*}
As before, the summation over $\bomega$ will vanish unless $t_{A'} :=
\sum_{(j,k) \in A'} (-1)^k (t^{(k)}_j - t^{(k)}_{j+1}) = 0$ for all
equivalence classes $A' \in A/\sim$, in which case the sum equals
$N^{|A/\sim|}$.  In particular, if $A/\sim$, the sum vanishes because
of the constraint $t^{(k)}_j \neq t^{(k)}_{j+1}$, so we may just as
well restrict the summation to those equivalence classes that contain
no singletons.  In particular we have
\begin{equation}\label{class-card}
|A/\sim| \leq K(n+1) = m.
\end{equation}
To summarize
\begin{eqnarray}
  \label{eq:expected3b}
\nonumber  
m(n,K) & =  & \sum_{\sim \in {\cal P}(A)} \,\, 
\sum_{\bt \in T^{2m}: t^{(k)}_j \neq t^{(k)}_{j+1} \text{ and } t_{A'} = 0
\text{ for all } A'} 
e^{i \sum_{k = 1}^{2K} (-1)^{k}
    \phi(t^{(k)}_{n+1})} 
N^{|A/\sim|} \prod_{A' \in A/\sim} F_{|A'|}(\tau)\\
& \le & \sum_{\sim \in {\cal P}(A)} \,\, 
\sum_{\bt \in T^{2K(n+1)}: t^{(k)}_j \neq t^{(k)}_{j+1} 
\text{ and } t_{A'} = 0
\text{ for all } A'} N^{|A/\sim|} \prod_{A' \in A/\sim} F_{|A'|}(\tau),  
\end{eqnarray}
since $|e^{i \sum_{k = 1}^{2K} (-1)^{k} \phi(t^{(k)}_{n+1})}| = 1$.
Observe the striking resemblance with \eqref{eq:expected3}.  Let
$\sim$ be an equivalence which does not contain any singleton.  
Then
the following inequality holds
$$
\# \,\, \{ \bt \in T^{2K(n+1)}: t_{A'} = 0, \text{ for all
} A' \in {A/\sim} \} \le |T|^{2K(n+1) - |A/\sim|}. 
$$
To see why this is true, observe as linear combinations of the
$t^{(k)}_j$ and of $t_0$, we see that the expressions $t^{(k)}_j -
t^{(k)}_{j+1}$ are all linearly independent, and hence the expressions
$\sum_{(j,k) \in A} (-1)^k (t^{(k)}_j - t^{(k)}_{j+1})$ are also
linearly independent.  Thus we have $|A/\sim|$ independent
constraints in the above sum, and so the number of $t$'s obeying the
constraints is bounded $|T|^{2n-|A/\sim|}$.

With the notations of section \ref{sec:moments}, we established
\begin{equation}
  \label{eq:expected52}
   m(n,K)   \le   \sum_{k = 1}^{m}  N^{k} \, |T|^{2m - k} \, P(2m,k) \, 
\sup_{\sim: |A/\sim| = k} \,\, \prod_{A' \in A/\sim} F_{|A'|}(\tau). 
\end{equation}
Now this is exactly the same as \eqref{eq:expected4} which we proved
obeys the desired bound.


\end{document}